\documentclass[a4paper,psamsfonts,reqno]{amsart}
\usepackage{amsmath}
\usepackage{amssymb} 
\usepackage[all]{xy} 
\usepackage{color}
\definecolor{Blue}{rgb}{0.3,0.3,0.9}
\definecolor{Red}{rgb}{0.9,0.3,0.3}

\begin{document}
\hyphenation{semi-per-fect Grothen-dieck}
\newtheorem{Lemma}{Lemma}[section]
\newtheorem{Th}[Lemma]{Theorem}
\newtheorem{Prop}[Lemma]{Proposition}
\newtheorem{OP}[Lemma]{Open Problem}
\newtheorem{Cor}[Lemma]{Corollary}
\newtheorem{Fact}[Lemma]{Fact}
\newtheorem{Def}[Lemma]{Definition}
\newtheorem{Ex}[Lemma]{Example}
\newtheorem{Exs}[Lemma]{Examples}
\newtheorem{Rem}[Lemma]{Remark}
\newenvironment{Remarks}{\noindent {\bf Remarks.}\ }{}
\newtheorem{Remark}[Lemma]{Remark}
\newenvironment{Proof}{{\sc Proof.}\ }{~\rule{1ex}{1ex}\vspace{0.5truecm}}
\newenvironment{ProofTh}{{\sc Proof of the Theorem.}\ }{~\rule{1ex}{1ex}\vspace{0.5truecm}}
\newcommand{\End}{\mbox{\rm End}}
\newcommand{\K}{\mbox{\rm K.dim}}
\newcommand{\Hom}{\mbox{\rm Hom}}
\newcommand{\infsupp}{\mbox{\rm inf-supp}\,}
\newcommand{\Ext}{\mbox{\rm Ext}}
\newcommand{\supp}{\mbox{\rm supp}\,}
\newcommand{\Max}{\mbox{\rm Max}}
\newcommand{\cl}{\mbox{\rm cl}}
\newcommand{\add}{\mbox{\rm add}}
\newcommand{\Inv}{\mbox{\rm Inv}}
\newcommand{\rk}{\mbox{\rm rk}}
\newcommand{\Tr}{\mbox{\rm Tr}}
\newcommand{\Sat}{\mbox{\bf Sat}}
\newcommand{\card}{\mbox{\rm card}}
\newcommand{\Ann}{\mbox{\rm Ann}}
\newcommand{\proj}{\mbox{\rm proj-}}
\newcommand{\codim}{\mbox{\rm codim}}
\newcommand{\B}{\mathcal{B}}
\newcommand{\Scal}{\mathcal{S}}
\newcommand{\Cong}{\mbox{\rm Cong}}
\newcommand{\Spec}{\mbox{\rm Spec-}}
\newcommand{\coker}{\mbox{\rm coker}}
\newcommand{\Cl}{\mbox{\rm Cl}}
\newcommand{\Ses}{\mbox{\rm Ses}}
\newcommand{\im}{\mbox{\rm Im}}
\newcommand{\Cal}[1]{{\mathcal #1}}
\newcommand{\+}{\oplus}
\newcommand{\N}{\mathbb N}
\newcommand{\No}{{\mathbb N}_0}
\newcommand{\Z}{\mathbb{Z}}
\newcommand{\Q}{\mathbb{Q}}
\newcommand{\C}{\mathbb{C}}
\newcommand{\T}{\mathbb{T}}
\newcommand{\R}{\mathbb{R}}
\newcommand{\notsim}{\sim\makebox[0\width][r]{$\slash\;\,$}}
\newcommand{\Mod}{\mbox{\rm Mod-}}
\newcommand{\lmod}{\mbox{\rm -mod}}
\newcommand{\mspec}{\mbox{\rm Max}}

\title[Reconstructing projective modules] {Reconstructing projective modules\\ from its trace ideal}

\author{Dolors Herbera}

\address{Departament de Matem\`atiques,
Universitat Aut\`onoma de Barcelona,  08193 Bellaterra
(Barcelona), Spain}
\email{dolors@mat.uab.cat}

\thanks{The first author wants to thank the hospitality of the Department of Algebra of  Charles University
(supported by GA\v CR P201/12/G028) where part of this paper was written. She was also 
partially supported by DGI MICIIN (Spain) through
Project MTM2011-28992-C02-01, and by the Comissionat per
Universitats i Recerca de la Generalitat de Catalunya through
Project 2005SGR00206}

  \author{Pavel P\v
r\'\i hoda}
\address{Charles University, Faculty of Mathematics and Physics \\Department
of Algebra, Sokolovsk\'a~83,
18675 Praha 8, Czech Republic}
\email{prihoda@karlin.mff.cuni.cz}

\thanks{The second author was supported by GA\v CR 201/09/0816  and research project GA\v CR P201/12/G028.}

\begin{abstract} We make a detailed study of idempotent ideals that are traces of countably generated projective right modules. We associate to such ideals an ascending chain of finitely generated left ideals and, dually, a descending chain of  cofinitely generated right ideals.  

The study of the first sequence allows us to characterize trace ideals of projective modules and to show that projective modules can always be lifted modulo the trace ideal of a projective module. As a consequence we give some new classification results of (countably generated) projective modules over  particular classes of semilocal rings. The study of the second sequence  leads us to consider projective modules over  noetherian FCR-algebras; we make some constructions of non-trivial projective modules  showing  that over such rings  the behavior of countably generated projective modules that are not direct sum of finitely generated ones is, in general, quite complex. 
\end{abstract}

\maketitle

It was proved by Whitehead \cite{whitehead} that an idempotent ideal that is finitely generated as a left ideal is the trace of a countably generated projective right $R$-module. For example, if $R$ is a left noetherian ring then any idempotent ideal is the trace ideal of a countably generated projective right $R$-module. 
P\v
r\'\i hoda in \cite{fair}   developed a machinery to work with countably generated projective modules in the setting of noetherian rings. These tools   have been extremely helpful in describing non-finitely generated projective modules over suitable classes of noetherian rings with low Krull dimension and not too many idempotent ideals.  See \cite{fair} for applications to integral group algebras and to some enveloping algebras of Lie algebras, also in \cite{PP, PP2} generalized Weyl algebras and some examples of lattices are studied. Finally, we mention \cite{HP} where the behavior of projective modules over semilocal noetherian rings is completely described.

One  of the main ideas in P\v r\'\i hoda's techniques is that  projective modules over suitable rings can be reconstructed  by determining  trace ideals and the finitely generated projective modules modulo such trace ideals.  
In this paper we extend   these results on trace ideals of projective modules  to a non necessarily noetherian  setting.  However, the not so exciting news are that  certainly further ideas will be needed to be able to have such complete classifications of projective modules also for some classes of non-noetherian rings. For example, just to understand completely projective modules over semilocal rings seems to be quite a hard problem. 

To explain which are the particular properties of trace ideals of projective modules, let us make some observations on trace ideals of the finitely generated ones. Let $P$ be a finitely generated projective right module over a ring $R$, and let  $E$ be an idempotent $n\times n$ matrix with entries in $R$ such that $P\cong ER^n$. The left ideal $J$ generated by the entries of $E$ is finitely generated and satisfies that $J^2=J$ and the same happens if we consider $K$ to be the right ideal of $R$ generated by the entries of $E$. Moreover, if $I$ is the trace of $P$ then $I=JR=RK$. 

Let  $I$ denote now the trace ideal of a countably generated  projective right $R$-module $P$. Let $E$ be a column-finite countable idempotent matrix defining $P$, then considering the left ideals generated by the entries in the first  columns of $E$ one constructs an ascending chain of finitely generated left ideals $J_1\subseteq J_2\subseteq \cdots \subseteq J_n\subseteq \cdots $ such that $J_{n+1}J_n=J_n$  and $I=\bigcup _{n\ge 1}J_nR=I$. The existence of such chains characterizes the ideals that are traces of countably generated projective right $R$-modules, cf. Proposition~\ref{traces}. In  Proposition~\ref{chartraces}, we show that trace ideals of arbitrary projective right $R$-modules can be characterized in terms of the existence of a direct system of ideals having such ascending sequences of finitely generated left ideals. As a consequence, we prove in Corollary~\ref{liftingproj}   that projective modules always lift modulo   trace ideals of   projective module. 

Further applications are developed in section~\ref{sectionsemisimperfect} where we describe countably generated projective modules over a particular class of semilocal rings that seems to be specially close to the class of semiperfect rings: the class of semilocal rings such that for any simple right  $R$-module $V$ there exists $n_V\ge 1$ such that $V^{n_V}$ has a projective cover. We call such a ring semi-semiperfect.  We give examples showing that, even being close to be semiperfect, semi-semiperfect rings can have countably generated projective modules that are not a direct sum of finitely generated ones; cf.~Example~\ref{ex:semisemiperfect}. 

Going back again to a   column-finite countable idempotent matrix $E$ defining a countably generated projective right $R$-module $P$,  we consider the right ideals generated by the entries in all rows of $E$, except perhaps   a finite number of them. With these ideals it is possible to construct a descending chain of right ideals $K_1\supseteq K_2\supseteq\cdots \supseteq K_n\supseteq \cdots $ such that $K_{n+1}K_{n}=K_{n+1}$ and $P/PK_n$ is finitely generated. The existence and the properties of  this descending chain gives  us a curious phenomena: if there exists an ideal $K$ of $R$ minimal with respect to the property $P/PK$ is finitely generated then there exists $n$ such that $K=RK_n$ and, hence, $K$ is an idempotent ideal. When this happens $P$ is said to be a fair sized projective module \cite{fair}. These ideas and some  consequences are developed in section~\ref{sectionfair}. Among other things we prove that  over  $R$ a semilocal ring  any projective module $P$ is fair-sized (Proposition~\ref{minimal}), but we fail to understand the r\^ole of the idempotent ideal $K$ associated to  $P$.

If $R$ is a noetherian ring (semilocal or not) and $P$ is a fair-sized projective module then the idempotent ideal $K$ is also the trace of a projective module. This can be quite helpful in determining the structure of $P$. For example, if $R$ is a noetherian ring  such that all countably generated projective right modules that  are not direct sums of finitely generated ones are fair sized, and there is only a finite number of idempotent ideals then, as it was shown in \cite{fair}, the projective modules over such rings are determined by pairs $(\overline{P}, I)$ where $I$ is an idempotent ideal (hence the trace of a projective module) and $\overline{P}$ is a finitely generated $R/I$-projective module. 

However,  over a noetherian ring not all countably generated projective modules that  are not direct sums of finitely generated ones are  fair-sized. In sections~\ref{sectionconstruction} and \ref{sectionfcr} we review the construction of a projective module with prescribed trace ideal due to Whitehead~\cite{whitehead}  and we prove that   if $R$ has an infinite descending chain of idempotent ideals with semisimple factors and starting at $I$, then $R$ has uncountably many non-isomorphic, non fair sized projective modules with  trace ideal $I$. Nice examples where this result can be applied are the enveloping algebras of  finite dimensional semisimple Lie algebras over an algebraic closed field and also  their quantum deformations. More generally,  noetherian $FCR$ algebras are the right context for these applications.   As recalled in Lemma~\ref{fglie}, any nonzero finitely generated projective module over the  enveloping algebra of  a finite dimensional semisimple Lie algebra is a generator and, hence, its trace is the whole ring. Therefore the projective modules in our construction are not direct sums of finitely generated ones, so it  looks really challenging to try to describe all projective modules over such   rings.

\section{On Whitehead's characterization of countably generated
projective modules}

All our rings are associative with $1$, and ring morphism means
unital ring morphism.

We think on (countably generated) projective modules as (countable) direct limits of finitely generated projective modules (i.e. as flat modules) that are  projective. To characterize when this happens it is important to keep in mind the following well known result. A complete proof of it can be found in \cite[Proposition
2.1]{angelerisaorin}.

\begin{Lemma}\label{split}  Let $R$ be a ring. Let
$M_1\stackrel{f_1}{\to} M_2\stackrel{f_2}{\to}\cdots
M_n\stackrel{f_n}{\to}M_{n+1}\cdots$ be a countable direct system
of right $R$-modules such that  $M=\varinjlim M_n$. Then the
canonical presentation of the direct limit
$0\to\oplus_{n\ge 1} M_n \to \oplus_{n\ge 1} M_n \to M\to 0$
splits if and only if there exists a sequence of homomorphisms
\[\cdots M_{n+1}\stackrel{g_n}{\to}M_n\to \cdots \to
M_2\stackrel{g_1}{\to}M_1\] such that:

For any $n\ge 1$ and for any finite subset $Y\subseteq M_n$ there
exists $\ell=\ell(n,Y)>n$ such that, for any $m> \ell$ and any $y\in
Y$, $g_m f_m\cdots f_n(y)=f_{m-1}\cdots f_n(y)$
\end{Lemma}

Our construction of countably generated projective modules relies on
Proposition \ref{directsystem} which is due to Whitehead
\cite[Theorem 1.9]{whitehead}. We want to explain how this result
can be seen as a consequence of the characterization of countably
generated projective modules as countably generated flat
Mittag-Leffler modules \cite{RG}.  To clarify the relation between
the two approaches we formulate the following lemma.

\begin{Lemma} \label{stationary} \emph{(Raynaud, Gruson \cite{RG})} Let $R$ be a ring. Let
$R^{m_1}\stackrel{f_1}\to R^{m_2} \stackrel{f_2}\to\cdots
R^{m_k}\stackrel{f_k}\to R^{m_{k+1}}\cdots $ be a countable direct
system of finitely generated free right modules with limit $P$. Then the following statements are equivalent:
\begin{itemize}
\item[(i)] For any $k\ge 1$ the descending chain of abelian groups
\[\mathrm{Hom}_R(R^{m_{k+1}},R)f_k\supseteq
\mathrm{Hom}_R(R^{m_{k+2}},R)f_{k+1}f_k\supseteq \cdots\supseteq
\mathrm{Hom}_R(R^{m_{k+\ell +1}},R)f_{k+\ell}\cdots f_k\supseteq
\cdots\] is stationary.
\item[(ii)] There exists a sequence of natural numbers $(\ell _k)_{k\ge 1}$ and a sequence of module
homomorphisms $(g_k\colon R^{m_{k+\ell _k+1}}\to R^{m_{k+\ell
_k}})_{k\ge 1}$ such that,  for any $k\ge 1$ and any $n\ge k+\ell
_k$, $g_nf_{n}\cdots f_k=f_{n-1}\cdots f_k$.
\item [(iii)] For any right module $M$, and for any $k\ge 1$ the descending chain of abelian groups
\[\mathrm{Hom}_R(R^{m_{k+1}},M)f_k\supseteq
\mathrm{Hom}_R(R^{m_{k+2}},M)f_{k+1}f_k\supseteq \cdots\supseteq
\mathrm{Hom}_R(R^{m_{k+\ell +1}},M)f_{k+\ell}\cdots f_k\supseteq
\cdots\] is stationary.
\item [(iv)] $P$ is
projective.
\end{itemize}
 Moreover, when the above equivalent statements hold, any direct system of finitely generated
projective modules with limit $P$ and of the form
\[F_1\stackrel{f'_1}\to F_2 \stackrel{f'_2}\to\cdots
F_k\stackrel{f'_k}\to F_{k+1}\cdots \] satisfies conditions $(i)$,
$(ii)$ and $(iii)$.
\end{Lemma}

\begin{Proof} Assume (i). Fix $k\ge 1$, and set $\ell _k\ge 1$
such that $\mathrm{Hom}_R(R^{m_{k+\ell _k +1}},R)f_{k+\ell _k}\cdots
f_k =\mathrm{Hom}_R(R^{m_{k+\ell _k }},R)f_{k+\ell _k-1}\cdots f_k$.

For any $i\in\{1, \dots, m_{k+\ell _k} \}$, let $\pi _i\colon
R^{m_{k+\ell _k }}\to R$ denote the projection onto the $i$-th
component, and let $\varepsilon _i \colon R\to R^{m_{k+\ell _k }}$
denote the canonical inclusion defined by $\varepsilon
_i(r)=(0,\dots, r^{i)},\dots ,0)$ for any $r\in R$. By hypothesis,
there exists $\omega _i\in \mathrm{Hom}_R(R^{m_{k+\ell _k +1}},R)$
such that
\[\pi _i f_{k+\ell _k-1}\cdots f_k=\omega _i f_{k+\ell _k}\cdots
f_k.\] Now $g_k=\sum _{i=1}^{m_{k+\ell _k }}\varepsilon _i\omega _i$
satisfies the desired properties. 

It is clear that $(ii)$ implies $(iii)$, and that $(iii)$ implies
$(i)$.  Using Proposition \ref{split} it is easy to see that  $(ii)$ holds if and only if $P$ is projective.

The final part of the statement  follows in a similar way
\cite{RG}.
\end{Proof}

\begin{Lemma}\label{strictml} Let $R$ be a ring. Let
$R^{m_1}\stackrel{f_1}\to R^{m_2} \stackrel{f_2}\to\cdots
R^{m_k}\stackrel{f_k}\to R^{m_{k+1}}\cdots $ be a countable direct
system of finitely generated free right modules with limit $P$. For any $k\ge 1$, let $u_k\colon R^{m_k}\to P$ denote the canonical
map. Then the following statements are equivalent:

\begin{itemize}
\item[(i)] For any $k>1$ there exist a  module homomorphism $g_k\colon R^{m_{k+1}}\to
R^{m_{k}}$ such that $g_kf_kf_{k-1}=f_{k-1}$.
\item[(ii)] For any $k\ge 1$,
$\mathrm{Hom}_R(P,R)u_k=\mathrm{Hom}_R(R^{m_{k+1}},R)f_k.$
\end{itemize}

When the above statements hold  $$\mathrm{Tr} (P)=\sum _{k\ge 1}
\left(\sum _{i=1}^{m_{k+1}}R\pi _if_k(R^{m_k})\right)$$ where $\pi
_i\colon R^{m_{k+1}}\to R$ denotes the projection onto the $i$-th
component. In particular if, for any $k\ge 1$, $f_k$ is given by
left multiplication by the matrix $X_k=(x_{ij}^k)$ then
$\mathrm{Tr} (P)=\sum _{k,i,j}Rx_{ij}^kR.$ 
\end{Lemma}

\begin{Proof} $(i)\Rightarrow (ii)$. We only need to prove that for any $k\ge 1$,
\[\mathrm{Hom}_R(P,R)u_k\supseteq\mathrm{Hom}_R(R^{m_{k+1}},R)f_k\]
Fix $k\ge 1$ and $\omega \in \mathrm{Hom}_R(R^{m_{k+1}},R)$. By
$(i)$, it follows that the sequence of homomorphisms $\{\omega
f_k,\omega g_{k+1}f_{k+1},\omega g_{k+1}g_{k+2}f_{k+2},\dots \}$
induces a homomorphism $f\colon \varinjlim R^{m_{k+\ell}}=P\to R$
such that $fu_k=\omega f_k$.

$(ii)\Rightarrow (i)$. Since for any $k\ge 1$, $\mathrm{Hom}_R(P,R)u_k\subseteq \mathrm{Hom}_R(R^{m_{k+2}},R)f_{k+1}f_k$, it follows from $(ii)$ that
$\mathrm{Hom}_R(R^{m_{k+1}},R)f_k=\mathrm{Hom}_R(R^{m_{k+2}},R)f_{k+1}f_k$.
 Therefore Lemma \ref{stationary}(i) holds and, for
any $k\ge 1$, the $\ell _k$ of Lemma \ref{stationary}(ii) is equal
to $1$. Hence $(i)$ holds.

It remains to prove the statement on the trace of $P$. It is clear that $\sum _{i=1}^{m_{k+1}}R\pi _if_k(R^{m_k})\subseteq
\mathrm{Tr} (P)$. Let $f\in\mathrm{Hom}_R(P,R)$. For any $p\in P$
there exists $k\ge 1$ and $x\in R^{m_k}$ such that $p=u_k(x)$.
Hence, by $(ii)$, there exists $\omega \in
\mathrm{Hom}_R(R^{m_{k+1}},R)$ such that
\[f(p)=f\circ u_k (x)=w\circ f_k(x)\in\sum _{i=1}^{m_{k+1}}R\pi
_if_k(R^{m_k}).\]
\end{Proof}

Now we  sketch the proof of Whitehead's
characterization.

\begin{Prop} \label{directsystem}\emph{\cite[Theorem 1.9]{whitehead}} Let $R$ be a ring, and let $P_R$ be a countably
presented (or generated) flat right  $R$-module. Then the following
statements are equivalent:
\begin{itemize}
\item[(i)] $P_R$ is projective.

\item[(ii)] There exists a direct system of finitely generated free
modules
\[R^{m_1}\stackrel{X_1}\to R^{m_2} \stackrel{X_2}\to\cdots
R^{m_k}\stackrel{X_k}\to R^{n_{k+1}}\cdots \] with limit $P$, and
where $X_k\colon R^{m_k} \to R^{m_{k+1}}$ denotes the homomorphism given by left multiplication by the matrix $X_k$, and a
sequence of matrices $\{Y_k\}_{k> 1}$ such that $Y_k\in
M_{m_{k}\times m_{k+1}}(R)$ and satisfies that $Y_kX_{k}X_{k-1}=X_{k-1}$ for
any $k> 1$.
\end{itemize}
In this situation if, for any $k\ge 1$, $X_k=(x^k_{ij})$  then
$\mathrm{Tr}\, (P)=\sum_{k\, , i\, , j} Rx^k_{ij}R$.
\end{Prop}

\begin{Proof} Assume $(i)$. Write $P$ as the direct limit of
\[F_1\stackrel{f_1}\to F_2 \stackrel{f_2}\to\cdots
F_k\stackrel{f_k}\to F_{k+1}\cdots \] 
where $F_k$ is finitely generated and free for any $k\ge 1$. By Lemma~\ref{stationary},  there exist a
sequence of natural numbers $(\ell _k)_{k\ge 1}$ and a sequence of
module homomorphisms $(g_k\colon F_{k+\ell _k+1}\to F_{k+\ell
_k})_{k\ge 1}$ such that,  for any $k\ge 1$ and any $n\ge k+\ell
_k$, $g_nf_{n}\cdots f_k=f_{n-1}\cdots f_k$. As done by Whitehead in
\cite[Theorem 1.9]{whitehead} (or also as in \cite{bazher}), this
allows to rearrange the direct system above to one with the same
limit and with the properties claimed in $(ii)$.

Lemma~\ref{stationary}, also shows that $(ii)$ implies $(i)$.

The rest of the statement follows from Lemma \ref{strictml}.
\end{Proof}

\section{An ascending chain of finitely generated left ideals}

Next lemma explains how to construct an ascending chain of finitely
generated left ideals out of a column finite idempotent matrix of
countable size.

\begin{Lemma}\label{remtrace} Let $R$ be a ring, and let $I$ be a
two-sided ideal of $R$ which is the trace ideal of a countably
generated projective right $R$-module.  Then there exists an
ascending chain of finitely generated left ideals $J_1\subseteq
J_2\subseteq \cdots \subseteq J_n\subseteq \cdots$ satisfying that,
for any $n\ge 1$, $J_{n+1}J_n=J_n$ and such that $I=\bigcup _{n\ge
1}J_nR$.
\end{Lemma}

\begin{Proof} Let $P$ be the countably generated
projective right module such that $I=\mathrm{Tr}\, (P)$. Let
$A=(a_{ij})_{i,j\ge 1}$ be a column finite countable idempotent
matrix with entries in $R$ such that $P=AR^{(\N)}$.  Note that
$I=\sum _{i,j\ge 1}Ra_{ij}R$.

For each $n\ge I$ consider the finitely generated left ideal
$L_n=\sum _{j\le n}Ra_{ij}$. Since $A$ is a column finite idempotent
matrix, for each $n\ge 1$, there exists $k_n> n$ such that
$L_{k_n}L_n=L_n$. Therefore by choosing a suitable  subchain of
$\{L_n\}_{n\ge 1}$  we can  construct an ascending chain of finitely
generated left ideals $\{J_n\}_{n\ge 1}$ such that $J_{n+1}J_n=J_n$
and $I=\bigcup_{n\ge 1}J_nR$. This finishes the proof of the lemma.
\end{Proof}

Reformulating a result by Whitehead \cite[Theorem~2.5]{whitehead},
we shall see that the existence of an ascending chain as in
Lemma~\ref{remtrace} characterizes the ideals that are traces of
countably generated projective right $R$-modules.

It is useful to keep in mind the following lemma, as  it explains some
modifications that can be made in the ascending chains appearing in
Lemma~\ref{remtrace}.

\begin{Lemma}\label{adding} Let $R$ be a ring. Let $J_1\subseteq J_2$ be finitely generated left ideals
of $R$ satisfying that $J_2J_1=J_1$. For $i=1,2$, fix $A_i$ a finite
set of generators of $J_i$.
\begin{itemize}
\item[(i)] Let $X$ be a finite subset of $R$ such that $1\in X$. For $i=1,2$, set $$J'_i=\sum _{r\in X\, a\in A_i
}  Rar.$$  Then $J'_1\subseteq J'_2$  and $J'_2J'_1=J'_1$. Moreover,
for $i=1,2$, $J_i\subseteq J_iR=J_i'R$.
\item[(ii)] There exists a finite   $X\subseteq R$ with $1\in
X$ such that $J_1$ is generated by $B\cdot A_1=\{b\cdot a\mid a\in
A_1\mbox{ and } b\in B\}$ where $$B=A_2\cdot X=\{ar\mid a \in A_2, \,
r\in X\}.$$ Moreover, if we set $J_2'= \sum _{b\in B }
Rb$ then $J'_2J_1=J_1$ and $J_2R=J_2'R$.
\end{itemize}
\end{Lemma}

\begin{Proof} Statement $(i)$ is immediate. To prove $(ii)$ observe
that for any $a\in A_1$
\[a=\sum _{a'\in A_1}y_{a'}^aa'\]
where $y_{a'}^a\in J_2R$ so that $y_{a'}^a$ is a sum of elements of
the form $sbr$ with $b\in A_2$ and $r$, $s\in R$. Let $X$ be a finite set of $R$ that
contains $1$ and all the elements $r$ in the expressions of the
elements $y_{a'}^a$.  Then $X$ satisfies the claimed properties.
\end{Proof}

\begin{Lemma} \label{independent} Let $R$ be a ring. Let $J_1$ and $J_2$ be left ideals
of $R$ such that there exist $a_1,\dots ,a_\ell \in J_1$ satisfying
that $J_1=\sum
_{i=1}^\ell J_2a_i$. Set $a=\begin{pmatrix}a_1\\ \vdots \\
a_\ell
\end{pmatrix}$. Fix $k\ge 1$ and let
\[A=\left( \begin{array}{ccc} a&\dots &0\\
\vdots & \ddots ^{k)} &\vdots \\ 0& \dots & a \end{array}\right)\in
M_{k\cdot \ell \times k}(J_1).\] Then, for any $m\ge 1$, $M_{m\times
k}(J_1)=M_{m\times k\cdot \ell}(J_2)A$.
\end{Lemma}

\begin{Proof} Clearly, $M_{m\times k\cdot
\ell}(J_2)A\subseteq M_{m\times k}(J_1)$. To prove the reverse
inclusion let $B=(b_{ij})\in M_{m\times k}(J_1)$. We may assume that
$B$ has only one nonzero entry, $b_{i_0\, j_0}$ say.

As $J_1=\sum _{i=1}^\ell J_2a_i$,  $b_{i_0\, j_0}=\sum _{i=1}^\ell
b_ia_i$ with $b_i\in J_2$ for $i=1,\dots ,\ell$. Denote the entries
of $A$ by $a_{ij}$ then, as any column of $A$ contains $a$, there
exists $i_1\ge 0$ such that
\[a=\begin{pmatrix}a_{i_1+1\, j_0}\\ \vdots \\ a_{i_1 +\ell \,
j_0}\end{pmatrix}\] and all other entries of the $j_0$ column of $A$
are zero.  Let $C=(c_{ij})\in M_{m\times k\cdot \ell}(J_2)$ be the
matrix such that all its entries are zero except possibly for
$c_{i_0\, i_1+1}=b_1,\dots ,c_{i_0\, i_1+\ell}=b_\ell$. Then $C\cdot
A=B$, which shows that $B\in M_{m\times k\cdot \ell}(J_2)A$.
\end{Proof}

\begin{Prop}\label{traces} Let $R$ be a ring, and let $I$ be a two-sided ideal of
$R$. Then the following statements are equivalent:
\begin{itemize}
\item[(i)] $I$ is the trace ideal  of a countably generated projective
right  $R$-module.
\item[(ii)] There exists an ascending chain of finitely generated
left ideals $\{J_n\}_{n\ge 1}$ such that $J_{n+1}J_n=J_n$ and
$I=\bigcup_{n\ge 1}J_nR$.
\item[(iii)] There exists an ascending chain of finitely generated
left ideals $\{J'_n\}_{n\ge 1}$ such that $I=\bigcup_{n\ge
1}J'_nR$, together with a suitable choice of a finite set of
generators $B_n$ of $J_n'$ such that $B_{n+1}\cdot B_n=\{b\cdot
b'\mid b\in B_{n+1}\mbox{ and } b'\in B_n\}$ generates $J'_n$. In
particular, $J'_{n+1}J'_n=J'_n$.
\end{itemize}
\end{Prop}

\begin{Proof} $(i)\Rightarrow (ii)$  and $(iii)\Rightarrow (i)$ follow from
\cite[Theorem~2.5]{whitehead}. Lemma~\ref{remtrace} gives a
different approach to $(i)\Rightarrow (ii)$. We also include a proof
of $(iii)\Rightarrow (i)$.

For every $n\ge 1$, set $Y_n=\begin{pmatrix}b^n_1\\ \vdots \\
b^n_{\ell_n}
\end{pmatrix}$. Define $X_1=Y_1$, and for $n>1$
\[X_n=\left( \begin{array}{ccc} Y_n&\dots &0\\
\vdots & \ddots{\ell _{n-1}\cdots \ell _1)} &\vdots \\ 0& \dots
&Y_n
\end{array}\right).\]
By Lemma \ref{independent}, for any $n\ge 1$
\[X_n\in M_{(\ell _n\cdots \ell _1)\times (\ell _{n}\cdots \ell
_1)}(J'_{n+1})X_n=M_{(\ell _{n}\cdots \ell _1)\times (\ell
_{n+1}\cdots \ell _1)}(R)X_{n+1}X_n.\] So that, for any $n\ge 1$,
there exists a matrix $C_n$ such that $C_nX_{n+1}X_n=X_n$.

By Proposition \ref{directsystem}, the limit of the direct system
\[R\stackrel{X_1}{\to}R^{\ell _1}\stackrel{X_2}{\to} R^{\ell _2\ell
_1}\cdots\] is a countably generated projective module with trace
$I$.

$(ii)\Rightarrow (iii)$. Let $\{J_n\}_{n\ge 1}$ be an
ascending chain of finitely generated left ideals with the
properties claimed in $(ii)$. Then $\{J'_n\}_{n\ge 1}$ and
$\{B_n\}$ are constructed by suitably enlarging the left ideals
$J_n$ by repeatedly applying Lemma~\ref{adding}.
\end{Proof}

Our  aim now is to give a characterization of trace ideals of arbitrary projective modules. Next lemma gives a more intrinsic approach to the construction of the sequence of Lemma~\ref{remtrace}.

\begin{Lemma} \label{finitelyidempotent} Let $R$ be a ring. Let $P$ be a projective right $R$-module and let $I=\mathrm{Tr}\, (P)$. Then, for any finite subset $X$ of $I$ there exist finitely generated left ideals $J_1\le J_2 \le I$ such that $X\subseteq J_1$ and $J_2J_1=J_1$.
\end{Lemma}

\begin{Proof} Since $X$ is finite, it is contained in $\mathrm{Tr}\, (P')$ where $P'$ is a countably generated direct summand of $P$. Thus, without lost of generality, we may assume that $P$ is countably generated. 

Fix $(p_i, \omega _i)_{i\ge 1}$ to be a dual basis of $P$. That is for any $i\ge 1$, $p_i\in P$ and $\omega _i\in \mathrm{Hom}_R(P,R)$, and for any $p\in P$, $p=\sum _{i\ge 1}p_i\omega _i(p)$ where $\omega_i(p)=0$ for almost all $i\ge 1$. Notice that, for any $f\in \mathrm{Hom}_R(P,R)$ and any $p\in P$, $f(p)=\sum _{i\ge 1}f(p_i)\omega _i(p).$ That is, there exists $n_p\ge 1$ such that the left ideal of $R$
\[J_p=\sum _{f\in \mathrm{Hom}_R(P,R)}f(p)=\sum _{i=1}^{n_p}R\omega _i(p)\le I\]
is finitely generated and satisfies that $J'_pJ_p=J_p$ where $J'_p=\sum _{i=1}^{n_p}J_{p_i}$; hence $J'_p$ is also finitely generated. 

Our previous argument shows that  $X\subseteq \sum _{p\in X} J_{p}=J_1$ and that  the left ideal $J_1$ is finitely generated; moreover  $J_2=J_1+\sum_{p\in X}  J'_{p}$  is also finitely generated, and $J_2J_1=J_1$.
\end{Proof}

\begin{Prop} \label{chartraces} Let $R$ be a ring, and let $I$ be an ideal of $R$. Then the following statement are equivalent:
\begin{itemize}
\item[(i)] There exists a projective right $R$-module $P$ such that $I=\mathrm{Tr}\, (P)$.
\item[(ii)] For any finite subset $X$ of $I$ there exists a couple of finitely generated left ideals $J_1\le J_2 \le I$ such that $X\subseteq J_1$ and $J_2J_1=J_1$.
\item[(iii)] There exists a left ideal $J\le I$, such that $JR=I$ and for any finite subset $X$ of $J$ there exist finitely generated left ideals $J_1\le J_2 \le J$ such that $X\subseteq J_1$ and $J_2J_1=J_1$.
\end{itemize}
Moreover, when the above equivalent statement hold, $I$ is an $\aleph _1$ directed union of ideals of the the form $\bigcup _{n\ge 1}J_nR$ where $J_1\le J_1\cdots \le J_n\le \cdots$ is a sequence of finitely generated left ideals such that $J_{n+1}J_n=J_n$ for any $n\ge 1$.
\end{Prop}

\begin{Proof} $(i)\Rightarrow (ii)$ follows from Lemma~\ref{finitelyidempotent}. The implication $(ii)\Rightarrow (iii)$ is trivial. Now we prove  $(iii)\Rightarrow (i)$.

Let $\Lambda=\{K\le J\mid K $ is finitely generated and there exists $K'\le J$  also finitely generated such that $K'K=K\}$.  Let $\mathcal{C}$ be the set of ascending chains in $\Lambda$ 
\[J_1\subseteq J_2\subseteq \cdots \subseteq J_n\subseteq \cdots \]
satisfying that $J_{n+1}J_n=J_n$ for any $n\ge 1$. By Proposition~\ref{traces}, for any $C\in \mathcal{C}$ there exists a countably generated projective right $R$-module $P_C$ such that $\mathrm{Tr}\, (P_C)=\sum _{J_n\in C}J_nR$.  Set $P=\sum _{C\in \mathcal{C}}P_C$, we claim that $I=\mathrm{Tr}\, (P)$. It is clear that $I\supseteq \mathrm{Tr}\, (P)$, to see the opposite inclusion, notice that, Lemma~\ref{finitelyidempotent}, for any finite subset $X$ of $J$ there exists $C\in \mathcal{C}$ such that $X\subseteq \sum _{J_n\in C}J_n$ and, hence $X\subseteq \mathrm{Tr}\, (P_C)$. This finishes the proof of the statement.

The proof of the remaining part of the statement follows from $(iii)\Rightarrow (i)$.
\end{Proof}

\begin{Cor} \label{tracefg} Let $R$ be a ring. Let $J$ be a finitely generated left
ideal such that $J^2=J$ then $JR$ is the trace of a (countably
generated) projective right $R$-module.
\end{Cor}

\begin{Proof} The ideal $I=JR$ fulfills
Proposition~\ref{traces}(ii) since the ascending chain can be
taken to be $J_n=J$ for any $n\ge 1$.
\end{Proof}

\begin{Cor} \label{noetherianjr} Let $R$ be a ring such that $R/J(R)$ satisfies the ascending chain condition on two-sided ideals
(e.g. $R/J(R)$ left or right noetherian). Then an ideal $I$ is the
trace of a (countably generated) projective right $R$-module if and
only if there exists a finitely generated left ideal $J$ such that
$J^2=J$ and $I=JR$.
\end{Cor}

\begin{Proof} Assume that $I$ is the trace ideal of  projective right $R$-module.
Let 
$\{J_n\}_{n\ge 1}$ be a sequence of finitely generated left ideals such that
$J_{n+1}J_n=J_n$ for any $n\ge 1$ and $I\supseteq \bigcup _{n\ge 1}J_nR$. We are going to see that such sequence is stationary.
Consider the ascending chain of two-sided ideals $\{J_nR\}_{n\ge
1}$. As $R/J(R)$ satisfies the ascending chain condition on
two-sided ideals, there exists $n_0$ such that
$J_{n_0}R+J(R)=J_{n_0+k}R+J(R)$ for any $k\ge 0$.

For any $n\ge n_0$,
\[J_n=(J_{n+1}+J(R))J_n=(J_{n+1}R+J(R))J_n=(J_{n}R+J(R))J_n=J_n^2+J(R)J_n.\]
As $J_n$ is finitely generated we can apply  Nakayama's Lemma to
deduce that  $J_n=J_n^2$ for any $n\ge n_0$.

Moreover, for any $n\ge n_0$,
\[J_n=J_n^2=(J_nR+J(R))J_n=(J_{n_0}R+J(R))J_n=J_{n_0}J_n+J(R)J_n\]
again Nakayama's Lemma allows us to deduce that $J_n=J_{n_0}J_n$.
Therefore, for any $n\ge n_0$, $J_nR=J_{n_0}J_nR=J_{n_0}R$.

 Now if $I$ is not of the form $JR$ for a finitely generated left ideal $J$ then Lemma~\ref{finitelyidempotent} would lead to the construction of a strictly ascending chain $\{J_n\}_{n\ge 1}$  of finitely generated left ideals such that
$J_{n+1}J_n=J_n$ and $J_nR\neq J_{n+1}R$ which is impossible by the above argument. Therefore there exists a finitely generated left ideal $J$ such that $JR=I$.  By Lemma~\ref{finitelyidempotent}, there exists a finitely generated left ideals $J_1 \le J'_2\le I$ such that $J\le J_1$ and $J'_2J_1=J_1$. Again by Lemma~\ref{finitelyidempotent}, there exists finitely generated ideals $J_2'\le J_2\le J_3\le I$ such that $J_3J_2=J_2$. This way we construct an ascending chain $J_1\le J_2\le J_3\le \cdots$ which, by the previous argument, it is stationary at a finitely generated ideal $J_{n_0}$ such that $J_{n_0} ^2=J_{n_0}$ and, since $J_1\le J_{n_0}$, also $I=J_{n_0}R$.
\end{Proof}

\begin{Cor} Let $R$ be a ring, and let $I$ and $I'$ be ideals of $R$. Assume that $I$ is the trace of
a  projective right $R$-module. If
$I+J(R)\subseteq I'+J(R)$ then $I\subseteq I'$.

In particular, if $I'$ is also the trace ideal of a   projective right $R$ module then $I+J(R)=I'+J(R)$ if and
only if $I=I'$.
\end{Cor}

\begin{Proof} Choose an ascending
chain of finitely generated left ideals $\{J_n\}_{n\ge 1}$
 satisfying that, for any $n$, $J_{n+1}J_n=J_n$ and such that $I\supseteq \bigcup _{n\ge
 1}J_nR$.

 Fix $n\ge 1$, then
 \[J_n=(J_{n+1}+J(R))J_n=(I'+J(R))J_n=I'J_n+J(R)J_n,\]
 by Nakayama's Lemma, $J_n=I'J_n\subseteq I'$ for any $n\ge 1$. Hence,  $\bigcup _{n\ge 1}J_nR\subseteq I'$. Since by Proposition~\ref{chartraces}, $I$ is a directed union of unions of such sequences, we deduce that $I'\subseteq I$.
\end{Proof}

\begin{Rem} \label{pureideals}
$(1)$ It was proved by Sakhaev \cite{sakhaev} that $I$ is a countably generated pure
right ideal of  $R$ if and only if $I$ can be generated by a
sequence $(a_n)_{n\in \N}$ of elements of $R$ satisfying that
$a_{n+1}a_n=a_n$ for any $n\in \N$. Since in this situation $R/I$ is
flat and countably presented, it follows from a result of Jensen
that $R/I$ has projective dimension at most one, so that $I_R$ is
projective. In particular, $RI=\mathrm{Tr}\, (I)$ is the trace of a
countably generated projective right $R$-module. Notice that, for
such $RI$, the sequence of ideals given by Proposition~\ref{traces}
can be taken to be $J_n=Ra_n$ for any $n\in \N$.

$(2)$  Let $I$ be a pure right ideal of $R$. Consider
$$\mathcal{S}=\{L\le I\mid L \mbox{ is countably generated and pure
in $R$}\}.$$ In view of  $(1)$ (or also because $I_R$ is a flat Mittag-Leffler module),  the set $\mathcal{S}$ is directed, closed under
countable unions and $I=\sum_{L\in \mathcal{S}} L$. Since any  $L\in
\mathcal{S}$ is a projective right ideal, we deduce that $RI$ is the
trace ideal of the projective module $P=\oplus_{L\in \mathcal{S}} L$. So that,  for any pure right ideal $I$ the two-sided ideal $RI$ is the trace of a projective right $R$-module.

$(3)$ In general, trace ideals of (finitely generated) projective modules are neither right nor left pure. For example, let $k$ be a field and let $R=\begin{pmatrix}k&k\\ 0&k\end{pmatrix}$. The artinian ring $R$ has two idempotents $e_1=\begin{pmatrix}1&0\\ 0&0\end{pmatrix}$ and $e_2=\begin{pmatrix}0&0\\ 0&1\end{pmatrix}$. The indecomposable projective right $R$-modules are $P_1=e_1R=\begin{pmatrix}k&k\\ 0&0\end{pmatrix}$ and $P_2=\begin{pmatrix}0&0\\ 0&k\end{pmatrix}$, and the indecomposable projective left $R$-modules are $Q_1=Re_1=\begin{pmatrix}k&0\\ 0&0\end{pmatrix}$ and $Q_2=Re_2=\begin{pmatrix}0&k\\ 0&k\end{pmatrix}$. Therefore $P_1=\mathrm{Tr} (P_1)=\mathrm{Tr} (Q_1)$ is pure as a right ideal but nor as a left ideal, while $Q_2=\mathrm{Tr} (P_2)=\mathrm{Tr} (Q_2)$ is pure as a left ideal but not as a right ideal.

Now in the ring $R\times R$ the two-sided ideal generated by the idempotent $(e_1,e_2)$ is a trace ideal of a projective module that is neither right nor left pure.

$(4)$ Not all trace ideals
of projective modules are obtained as in $(2)$.
For example, Sakhaev constructed a ring $R$ such that $R_R/J(R)$ is a direct sum of $2$ non-isomorphic simple modules $S_1,S_2$, there is no projective 
$R$-module $Q$ such that $Q/QJ(R)$ is simple and there exists a projective module $P$ such that $P/PJ(R) \simeq S_1^2$; see also \cite{HP2}. We claim that such ring has no-nontrivial pure right ideals, and hence $\mathrm{Tr}\, (Q)\neq RI$ for a pure ideal $I$ of $R$. 

Indeed, if $R$ has a pure right ideal different from $R$ and $\{0\}$ then it has a non trivial countably generated    pure right ideal; call such an ideal $I$. By $(1)$, $I$ is projective.  Since $I/IJ(R)$ is a submodule of $R/J(R)$ and it is neither isomorphic to $S_1$ nor $S_2$, then 
  $I = 0$ or $I = R$, a contradiction.  \end{Rem}    
 
For a commutative ring $R$, Vasconcelos showed that  trace
ideals of projective modules over $R$ are precisely  the pure ideals. We give an alternative proof of Vasconcelos result.

\begin{Lemma}\label{determinant} Let $R$ be a commutative ring. Let $J_1\subseteq J_2$
be ideals such that $J_1$ is finitely generated and $J_2J_1=J_1$.
Then there exists $a\in J_2$ such that $ab=b$ for any $b\in J_1$.
\end{Lemma}

\begin{Proof} Let $b_1,\dots ,b_r$ be a finite set of generators of
$J_1$. By hypothesis, for any $i=1,\dots ,r$, there exists
$a_{ij}\in J_2$ with $i=1,\dots ,r$, such that $b_i=\sum
_{j=1}^ra_{ij}b_j$. So that, if  $A=(a_{ij})\in M_r(J_2)$ then
\[(\mathrm{Id_n}-A)\left(\begin{array}{c}b_1\\\vdots \\
b_r\end{array}\right)=0.\] Since
$\mathrm{det}\,(\mathrm{Id_n}-A)=1-a$ with $a\in J_2$, we deduce
that $ab_i=b_i$ for any $i=1,\dots ,r$. So that $a$ has the claimed
property.
\end{Proof}

\begin{Prop}\label{commcg} Let $R$ be a commutative ring, and let $I$ be an ideal of $R$. The following are
equivalent statements:
\begin{itemize}
\item[(i)] $I$ is the trace ideal of a countably generated
projective module.
\item[(ii)] $I$ can be generated by a sequence
$(a_n)_{n\in \N}$ of elements of $R$ satisfying that
$a_{n+1}a_n=a_n$ for any $n\in \N$.
\item[(iii)] $I$ is a countably generated pure ideal of $R$.
\item[(iv)] $I$ is  a countably generated, projective and   pure ideal
of $R$.
\end{itemize}
\end{Prop}

\begin{Proof} In view of Remark~\ref{pureideals}, it is enough to show
that $(i)\Rightarrow (ii)$. By Proposition~\ref{traces}, $I=\bigcup
_{n\in \N}J_n$ where $J_n$ are finitely generated ideals of $R$
satisfying that $J_{n+1}J_n=J_n$. By repeatedly using
Lemma~\ref{determinant} we find that, for any $n\in \N$, there exists $a_n\in
J_{n+1}$ such that $a_nb=b$ for any $b\in J_n$. The sequence
$(a_n)_{n\in \N}$ satisfies the desired properties.
\end{Proof}

\begin{Cor} Let $R$ be a commutative ring. An ideal $I$ of $R$ is
the trace of a projective $R$-module if and only if it is pure in
$R$.
\end{Cor}

\begin{Proof} In view of Remark~\ref{pureideals}, it is enough to show that if $I$ is an ideal of $R$ that is the trace
of a projective module $P$ then it is pure. Since any projective module is a direct
sum of countably generated projective submodules, we can write
$P=\bigcup _{\alpha \in \Lambda}P_\alpha$ where $(P_\alpha)_{\alpha
\in \Lambda}$ for a directed set of countably generated direct
summands of $P$. By Proposition~\ref{commcg}, $I=\mathrm{Tr}
(P)=\sum _{\alpha \in \Lambda}\mathrm{Tr} (P_\alpha)$ is a directed
union of pure ideals of $R$, so that it is a pure ideal of $R$.
\end{Proof}

\begin{Cor} Let $R$ be a commutative ring such that all pure ideals
of $R/J(R)$ are finitely generated (hence generated by an
idempotent). Then all pure ideals of $R$ are of the form $Re$ for
$e=e^2\in R$.
\end{Cor}

\begin{Proof} Let $I$ be a countably generated pure ideal of $R$.
By Remark~\ref{pureideals}, $I=\sum _{n\in \N}a_nR$ with
$a_{n+1}a_n=a_n$ for any $n\in \N$. By hypothesis, there exists
$n_0\in \N$, such that $a_{n_0}R+J(R)=a_{n_0+k}R+J(R)$ for any $k\in
\N$. By Corollary~\ref{noetherianjr}, $I=a_{n_0}R$. Since $R/I$ is a finitely presented flat module (or using the
determinant trick) we deduce that $I$ and $R/I$ are projective, so
that, $I$ is generated by an idempotent.

In particular,  an ascending chain
$$e_1R\subseteq e_2R\subseteq\cdots e_nR\subseteq\cdots $$
where, for any $n\ge 1$, $e_n^2=e_n\in R$ must be stationary because $I=\bigcup _{n\ge 1}e_nR$ is a pure ideal of $R$.

If $I$ is any pure ideal of $R$ then, as pointed out in
Remark~\ref{pureideals}, the set of countably generated pure ideals
contained in $I$ form a direct system, closed under countable unions
and with limit $I$. By the first part of the proof, $I$ is
generated by idempotents. Since any ascending chain of principal ideals generated by an idempotent must terminate, we deduce that $I$ is generated by a single idempotent.
\end{Proof}

\section{Lifting projective modules modulo a trace ideal}

Let $I$ be an ideal of $R$ that is the trace of a  projective right module. Let $J+I$ be an ideal of $R/I$
that is the trace ideal of a projective right
$R/I$-module. Using the characterization in
 Proposition~\ref{chartraces}, it is not difficult to see that $I+J$
 is also a trace ideal of a   projective right module
 over $R$. We will see the  stronger result that a projective
 module modulo a trace ideal can be lifted to a projective module
 over $R$.  The main step will be to prove the result in the countably generated case,  to do that we extend the ideas of \cite[Lemma~2.6]{fair} which,
 in turn, are based in Whitehead's method
 to construct countably
generated projective modules with  prescribed trace ideal; cf.  \cite{whitehead} or Lemma~\ref{strictml}.

\begin{Th} \label{liftingcountable} Let $R$ be a ring, and let $I$ be an ideal of $R$ that
is the trace of a countably generated projective right $R$-module.
Let $P'$ be a countably generated projective right module over
$R/I$, then there exists a countably generated projective right
$R$-module $P$ such that $I\subseteq \mathrm{Tr}\,
(P)$,   $P/PI\cong P'$ and $\mathrm{Tr}\,
(P)/I= \mathrm{Tr}_{R/I}
(P')$.\end{Th}

\begin{Proof} According to Proposition~\ref{directsystem}, we can choose $\{X_k\}_{k\ge 1}$ and $\{Y_k\}_{k\ge 1}$  two
sequences of finite matrices with entries in $R$ such that they
define a direct system of finitely generated free $R/I$-modules
\[(R/I)^{n_1}\stackrel{\overline{X_1}}\to (R/I)^{n_2} \stackrel{\overline{X_2}}\to\cdots
(R/I)^{n_k}\stackrel{\overline{X_k}}\to (R/I)^{n_{k+1}}\cdots \]
with limit $P'$ and such that $Y_kX_{k+1}X_k-X_k \in
M_{n_{k+1}\times n_k}(I)$ for any $k\ge 1$. To avoid trivial
situations we assume $n_k\ge 1$ for any $k\ge 1$.

By Proposition~\ref{traces}(iii) and Lemma~\ref{adding}, $I=\bigcup
_{k\ge 1}J_kR$ where $\{J_k\}_{k\ge 1}$ is an ascending chain of
finitely generated left ideals of $R$ such that, for each $k\ge 1$,
$Y_kX_{k+1}X_k-X_k \in M_{n_{k+1}\times n_k}(J_{k})$ and there
exists a nonempty finite set $\{a_1^k,\dots ,a_{\ell _k}^k\}$ of
generators of $J_k$ such that $J_k=\sum _{i=1}^{\ell
_k}J_{k+1}a^k_\ell$. For any $k\ge 1$, let $a_k=\begin{pmatrix}a^k_1\\ \vdots \\
a^k_{\ell_k} \end{pmatrix}$.

We define inductively a sequence of natural numbers
$\{m_k\}_{k\ge1}$ by setting  $m_1=n_1$ and, for $k\ge 1$,
$m_{k+1}=n_{k+1}+m_k\ell _k$. Now we want to construct a direct
limit of finitely generated free modules
\[R^{m_1}\stackrel{A_1}\to R^{m_2} \stackrel{A_2}\to\cdots
R^{m_k}\stackrel{A_k}\to R^{m_{k+1}}\cdots \] such that $\varinjlim
(R^{m_k},A_k)=P$ is a projective module and $P\otimes R/I\cong P'$.
To this aim and for each $k\ge 1$ we name the following auxiliary
matrices,
\[A'_k=\left(\begin{array}{ccc}a_k&\dots &0\\
&\ddots \!^{m_k)}& \\
0&\dots &a_k\end{array}\right)\in M_{(m_k \ell _k)\times
m_k}(R)\mbox{ and } X'_k=\left(\begin{array}{cc}X_k&
0\end{array}\right)\in M_{n_{k+1}\times m_k}(R).\]

Now we  set,
for each $k\ge 1$, $A_k=\left(\begin{array}{c}X'_k\\
A'_k\end{array}\right)\in M_{m_{k+1}\times m_k}(R)$. Observe that
\[\varinjlim (R^{m_k}, A_k)\otimes _R R/I\cong \varinjlim
\left( (R/I)^{m_k},\left(\smallmatrix \overline{X'_k}\\
0\endsmallmatrix \right)\right)\cong P'.\]

By Proposition~\ref{directsystem}, to conclude that $P$ is a
projective module, we need to find a sequence of square matrices
$\{C_k\}_{k\ge 1}$ such that $C_k\in M_{m_{(k+1)}\times
m_{(k+2)}}(R)A_{k+1}$ and $C_kA_k=A_k$ for any $k\ge 1.$

 Since, by construction,
$Z_k=Y_kX_{k+1}X'_k-X'_k \in M_{n_{k+1}\times m_k}(J_{k})$, by
Lemma~\ref{independent}, $Z_k=C_1^kA'_k$ for some $C_1^k\in
M_{n_{k+1}\times m_k\ell_k}(J_{k+1})$. Also by Lemma~\ref{independent},
there exists $C_2^k\in M_{m_k\ell _k}(J_{k+1})$ such that
$C_2^kA'_k=A'_k$. Then, for any $k\ge 1$, the matrix
\[C_k=\left(\begin{array}{cc}Y_kX_{k+1}&-C^k_1\\ 0&C^k_2\end{array}
\right)\] satisfies that $C_kA_k=A_k$. 

Write $A_{k+1} = \left(\begin{array}{cc} X_{k+1} & 0 \\ B_{k+1} & B_{k+1}' \end{array} \right)$.
Let $D_2^k \in M_{m_k\ell_k\times m_{k+1}\ell_{k+1}}(R)$ be such that $(0\ C_2^k) = D_2^kA_{k+1}'$. Further, 
using Lemma~\ref{independent} we find $D_1^k \in M_{n_{k+1}\times m_{k+1}\ell_{k+1}}(R)$ such that 
$D_1^kB_{k+1} = 0$ and $D_1^k B'_{k+1} = -C_1^k$. Then 
\[C_k = \left(\begin{array}{cc} Y_k & D_1^k \\ 0 & D_2^k \end{array} \right)A_{k+1}.\]
Hence $P$ is a projective module as claimed.

Finally, by Proposition~\ref{directsystem} and the
construction of $P$, it also follows that  $\mathrm{Tr}\,
(P)/I= \mathrm{Tr}_{R/I}
(P')$.
\end{Proof}

\begin{Cor} \label{liftingproj} Let $R$ be a ring, and let $I$ be an ideal of $R$ that
is the trace of a  projective right $R$-module.
Let $P'$ be a  projective right module over
$R/I$, then there exists a  projective right
$R$-module $P$ such that $I\subseteq \mathrm{Tr}\,
(P)$,   $P/PI\cong P'$ and $\mathrm{Tr}\,
(P)/I= \mathrm{Tr}_{R/I}
(P')$.\end{Cor}

\begin{Proof} We may assume that $P'$ is a countably generated projective right $R/I$-module. By Proposition~\ref{directsystem}, we can choose $\{X_k\}_{k\ge 1}$ and $\{Y_k\}_{k\ge 1}$  two
sequences of finite matrices with entries in $R$ such that they
define a direct system of finitely generated free $R/I$-modules
\[(R/I)^{n_1}\stackrel{\overline{X_1}}\to (R/I)^{n_2} \stackrel{\overline{X_2}}\to\cdots
(R/I)^{n_k}\stackrel{\overline{X_k}}\to (R/I)^{n_{k+1}}\cdots \]
with limit $P'$ and such that $Y_kX_{k+1}X_k-X_k \in
M_{n_{k+1}\times n_k}(I)$ for any $k\ge 1$.

By Proposition \ref{chartraces}, there exists $I_C\subseteq I$ such that $I_C$ is the trace of a countably generated projective right $R$-module, and $I_C$ contains the entries of the matrices $(Y_kX_{k+1}X_k-X_k)_{k\ge 1}$. Consider the directed system of free $R/I_C$-modules,
\[(R/I_C)^{n_1}\stackrel{\overline{X_1}}\to (R/I_C)^{n_2} \stackrel{\overline{X_2}}\to\cdots
(R/I)^{n_k}\stackrel{\overline{X_k}}\to (R/I_C)^{n_{k+1}}\cdots \]
By  Proposition~\ref{directsystem}, its direct limit is a countably generated projective right $R/I_C$-module, $P''$ say. By Theorem~\ref{liftingcountable}, there exists $P$ a countably generated projective right $R$-module such that $\mathrm{Tr}\, (P)\supseteq I_C$ and $P/PI_C\cong P''$.

Let $Q$ be a projective right $R$-module such that $\mathrm{Tr}\, (Q)=I$. Then $P\oplus Q$ satisfies the desired properties.
\end{Proof}

\section{An application: projective modules over semi-semiperfect rings} \label{sectionsemisimperfect}

A ring $R$ is said to  be semilocal if $R/J(R)$ is semisimple artinian. It is said to be semi-semiperfect if it is semilocal and  for any simple right $R$-module $V$ there exists $n_V\in \N$ and a finitely generated projective right $R$-module $P$ such that $P/PJ(R)\cong V^{n_V}$. Observe that, if for any simple right $R$-module $V$, $n_V$ can be taken to be $1$ then $R$ is a semiperfect ring. 

Semi-semiperfect rings are closed by quotients by two-sided ideals. For further quoting we record this fact in the next lemma.

\begin{Lemma} \label{semiquotient}  Let $R$ be a semilocal ring, and let $I\neq R$ be a two-sided ideal of $R$. If $R$ is semi-semiperfect then so is $R/I$.
\end{Lemma}

Let $R$ be a ring. We denote by $V(R)$ the (commutative) monoid of isomorphism classes of finitely generated projective right $R$-modules with the addition induced by the direct sum of projective modules.  We denote by $V^*(R)=V^*(R_R)$ the (commutative) monoid of isomorphism classes of countably generated projective right $R$-modules with the addition induced by the direct sum of projective modules.

We note the following corollary of Theorem~\ref{liftingcountable}.

\begin{Cor} \label{ontovstar} Let $R$ be a ring and let $I$ be the trace ideal of a countably generated projective module. Then the assignment $\langle P\rangle \mapsto \langle P/PI\rangle$ for $\langle P\rangle\in V^*(R)$ defines an onto morphisms of monoids $V^*(R)\to V^*(R/I)$.
\end{Cor}

As an application of the results in the previous sections we will show that if $R$ is a semi-semiperfect ring then the monoids $V^*(R)$ are finitely generated and, in fact, are included in the class of monoids that can be realized as $V^*(R)$ for $R$ a semilocal noetherian ring.

We introduce some notation in order to be able to state and prove our results. 

Our convention is $\N =\{1,2,\dots \}$, and we denote the
nonnegative integers by $\No =\{0, 1,2, \dots \}$. We also need the monoid
$(\mathbb{N}_0^*,+,0)$ whose underlying set is $\mathbb{N}_0 \cup
\{\infty\}$, the operation $+$ is the extension of addition of
non-negative integers by the rule $\infty+x = x + \infty = \infty$.
It is  also interesting to extend the product of
$\No$ to $\No^*$ by setting $\infty\cdot 0=0$ and $\infty \cdot
n=\infty$ for any $n\in \No^*\setminus \{0\}.$

Let  $R$ be a  semilocal ring such that $R/J(R)\cong
M_{n_1}(D_1)\times \cdots \times M_{n_k}(D_k)$ for suitable
 division rings $D_1,\dots ,D_k$. Fix an onto ring homomorphism
$\varphi \colon R\to  M_{n_1}(D_1)\times \cdots \times M_{n_k}(D_k)$
such that $\mathrm{Ker}\, \varphi =J(R)$.  

In this situation, the monoids $V(R)$ and $V^*(R)$ can be seen a submonoids of $\No^k$ and of 
 $(\No^*)^k$, respectively. We briefly explain how this can be made. We will keep using the fact that two projective right $R$-modules are isomorphic if and only if they are isomorphic modulo $J(R)$ \cite{P1}.

Fix $V_1,\dots ,V_k$ to be an ordered set of representatives of the isomorphism classes of simple right $R$-modules. For $i=1,\dots ,k$, we assume that $\mathrm{End}_R(V_i)=D_i$.  For any countably generated projective right
$R$-module $P$, $P/PJ(R)\cong V_1^{(I_1)}\oplus \cdots \oplus V_k^{(I_k)}$. Denote by $\langle P \rangle$ the isomorphism class of $P$, then we set 
\[\dim _\varphi (\langle P\rangle): =(m_1,\dots ,m_k)\in (\No^*)^k,\] 
where, for $i=1,\dots ,k$, $m_i=\vert I_i\vert$ if $I_i$ is finite and $m_i=\infty$ if $I_i$ is infinite. Since $P\cong Q$ if and only if $P/PJ(R)\cong Q/QJ(R)$ \cite{P1}, the map $\dim _\varphi\colon V^*(R)\to (\No^*)^k$ is an injective morphism of additive monoids. 

Notice that  $\dim _\varphi (V(R))\subseteq \No^k$. Moreover, in this case the monoid morphism $\dim _\varphi$ is an embedding of monoids (that is, if in $\No^k$ there is an equality $x+y=z$ with $x$, $z\in \dim _\varphi (V(R))$ then also $y\in \dim _\varphi (V(R))$). It was proved in \cite{FH} that this property characterizes the monoids that can be realized as $V(R)$ for a semilocal ring $R$.

The following result characterizes  semilocal rings that are semi-semiperfect in terms of the values of $\dim _\varphi$ and also shows that to test that a projective module is finitely generated it is enough to show that it is finitely generated modulo the Jacobson radical.

\begin{Lemma} \label{finitesupport} With the  notation above. The ring $R$ is semi-semiperfect if and only if, for $i=1,\dots ,k$, there exists $m_i\in \N$ such that $(0,\dots, m_i^{i)},\dots ,0)\in \dim _\varphi (V(R))$. 

Moreover, in this situation, if $P_R$ is a countably generated projective right $R$-module such that $P/PJ(R)$ is finitely generated then $P_R$ is finitely generated.
\end{Lemma}

\begin{Proof} The first part of the statement is just a translation of the definition of semi-semiperfect ring to the monoid $V(R)$.

To prove the second part, let $P_R$ be a countably generated projective right $R$-module such that $P/PJ(R)$ is finitely generated.  Hence there exists   $(a_1,\dots ,a_k)\in  \No^k$  such that $\dim _\varphi(\langle P\rangle )=(a_1,\dots ,a_k)$.  For $i=1,\dots ,k$, let $P_i$ be a finitely generated right projective module such that $\dim _\varphi (\langle P_i\rangle)=(0,\dots, m_i^{i)},\dots ,0)$ and set $t_i=\prod _{j\neq i}m_j$. Therefore  $$\dim _\varphi(\langle P^{m_1\cdots m_k}\rangle)=\dim _\varphi(\langle \oplus _{i=1}^kP_i^{a_it_i}\rangle)$$
so that $P^{m_1\cdots m_k}\cong  \oplus _{i=1}^kP_i^{a_it_i}$ is finitely generated and, hence, we conclude that  $P$ is a finitely generated projective right $R$-module.
\end{Proof}

Let $\mathbf{x}=(x_1,\dots ,x_k)\in (\No^*) ^k$. We define
\[\mathrm{supp}\, (\mathbf{x})=\{i\in \{1,\dots, k\}\mid x_i\neq 0\}\]
and we refer to this set as the \emph{support of $\mathbf{x}$}. We
also define \[\infsupp  (\mathbf{x})=\{i\in \{1,\dots, k\}\mid x_i=
\infty \},\]  we refer to this set as the \emph{infinite support of
$\mathbf{x}$}.

In the next lemma we recall that if $P$ is a projective right $R$-module then the trace ideal records the support of $\dim _\varphi(\langle P\rangle )$.

\begin{Lemma} \label{simplefactors} \emph{(\cite[Lemma~2.2]{HP})} Let $P$ be a  countably generated projective
right $R$-module with trace ideal $I$. Set $\mathbf{x}=\dim _{\varphi} (\langle P\rangle )$. For $i=1,\dots, k$, the following statement are equivalent:
\begin{itemize}
\item[(i)] $V_i$ is a quotient of $P$.
\item[(ii)] $V_i$ is a quotient of $I$.
\item[(iii)] $I+r_R(V_i)=R$.
\item[(iv)] $i\in \mathrm{supp}\, (\mathbf{x})$.
\end{itemize}
\end{Lemma}

In view of Lemma~\ref{simplefactors} there is an injective map between
\[\Phi\colon \mathcal{T}=\{I\le R\mid I\mbox{ is a trace ideal of a  projective right module}\}\to \mathcal{P}(\{1,\dots ,k\})\]
given by $\Phi (I)= \mathrm{supp}\, \dim _{\varphi} (\langle P\rangle )$ where $P_R$ is any (countably) generated projective right $R$-module with trace ideal $I$. Notice that if $R$ is semi-semiperfect then $\Phi$ is a bijective map.

\begin{Th} \label{semisemiperfect} Let $R$ be a semi-semiperfect ring.  For $i=1,\dots ,k$, let $P_i$ be a finitely generated projective right $R$-module such that $\Phi (\mathrm{Tr}\, (P_i))=\{i\}$. For any $J\in \mathcal{T}$   fix $G_J$ a finite set of projective right $R$-modules such that $\overline{G_J}=
\{Q/QJ\mid Q\in G_J\}$ generates $V(R/J)$.  

Let $P$ be a countably generated projective right $R$-module,  let $X=\infsupp \,( \dim _\varphi(\langle P\rangle ))$ and let $I$ be a trace ideal such that $\Phi (I)=X$ then
\[P\cong \left( \oplus _{i\in X}P_i^{(\omega _0)}\right)\oplus \left( \oplus _{Q\in G_I} Q^{n_Q}\right)\]
where, for any $Q\in G_I$, $n_Q\in \No$.

In particular, $V^*(R)$ is finitely  generated by $P_1^{(\omega _0)}, \dots , P_k^{(\omega _0)}$ and $\bigcup _{I\in \mathcal{T}} G_I$.
\end{Th}

\begin{Proof} The statement is trivial if $I=R$, so we may assume that $I\neq R$. By Lemma~\ref{semiquotient}, $R/I$ is a semi-semiperfect ring.  Let $\dim _{\varphi _I}$ be the dimension function induced by $\dim _\varphi$ over $R/I$. Such dimension function satisfies that for any countably generated projective right $R$-module $K$
 $$\dim _{\varphi _I} (\langle K/KI\rangle)=\pi _X(\dim _\varphi \, (\langle K\rangle)
 $$
 Where $\pi _X\colon (\No^*)^k\to (\No^*)^{\{1,\dots ,k\}\setminus X}$ denotes the canonical projection.

By Lemma~\ref{simplefactors} and because of  the definition of $I$, $\dim _{\varphi _I} (\langle P/PI\rangle)$ has all its components in $\N_0$ so that, by Lemma~\ref{finitesupport}, $P/PI$ is a finitely generated projective right $R/I$ module. Therefore  $P/PI\cong  \oplus _{Q\in G_I} (Q/QI)^{n_Q} $ for suitable $n_Q\in \No$.  Now it follows that 
\[P\cong K= \left( \oplus _{i\in X}P_i^{(\omega)}\right)\oplus \left( \oplus _{Q\in G_I} Q^{n_Q}\right)\]
because $\dim _\varphi (\langle P\rangle)= \dim _\varphi (\langle K\rangle)$, so that $P$ and $K$ are isomorphic modulo $J(R)$, hence, they are isomorphic.
\end{Proof}

\begin{Remark} Notice that the crucial fact in Theorem~\ref{semisemiperfect} is that any infinite support of a countably generated projective $R$-module is the support of a projective module.   A   result in the same spirit of Theorem~\ref{semisemiperfect} could be obtained with a semilocal ring $R$ such that for any simple right $R$-module $V$ there exists a projective right module $P_V$ such that $P_V/P_VJ(R)\cong V^{(I)}$ for a suitable, possibly infinite, non-empty set $I$. If $R/J(R)\cong M_{n_1}(D_1)\times \cdots \times M_{n_k}(D_k)$ for suitable
 division rings $D_1,\dots ,D_k$, then we would also have that the map $\Phi\colon \mathcal{T}\to  \mathcal{P}(\{1,\dots ,k\})$ is onto. However instead of lifting to $ V^*(R)$ the generators of $V(R/I)$ for any $I\in \mathcal{T}$ we would need to lift the generators of $W(R/I)$ to obtain a set of generators of the monoid $V^*(R)$.
 
 We recall that for a semilocal ring $R$, $W(R)$ is the submonoid of $V^*(R)$ consisting of the isomorphism classes of countably generated projective right $R$-modules that are finitely generated modulo $J(R)$. As stated in Lemma~\ref{finitesupport}, for a semi-semiperfect ring $W(R)=V(R)$ but this is not true in the more general setting that we are describing, see \cite[Example~3.6]{HP2}.
\end{Remark}

It is not difficult to see that the monoids described in Theorem~\ref{semisemiperfect} are a particular instance of the monoids given by a system of supports in the sense of \cite[\S 7]{HP}. Therefore, they can be described as the solutions in $\No^*$ of suitable systems of linear equations and linear congruences with coefficients in $\No$ \cite[Theorem~7.7]{HP}. Now we use this technology to give some examples that illustrate that for semi-semiperfect rings not all projective modules are a direct sum of finitely generated ones. 

\begin{Ex} \label{ex:semisemiperfect} Let $R$ be a semi-semiperfect ring and we fix an onto ring homomorphism $\varphi\colon R\to
 M_{n_1}(D_1)\times \cdots \times M_{n_k}(D_k)$ with kernel $J(R)$ and such that $D_1,\dots ,D_k$ are division rings. Then, by \cite[Proposition~6.2]{HP}, there   there exist $0\le n \le k$,  $A\in M_{n\times k}(\No)$ and $m_1,\dots ,m_n$
integers strictly bigger than one such that such that $\mathbf{x}\in \dim _\varphi (V(R))$ if and only if 
\begin{equation}  A\cdot \mathbf{x}^t\in \left(\begin{array}{c} m_1\No\\ \vdots \\ m_n\No \end{array}\right) \end{equation}

By \cite[Theorem~2.6]{HP}, $R$ can be chosen such that it is noetherian semilocal and $\mathbf{x}\in \dim _\varphi (V^*(R))$
if and only if 
\begin{equation}  A\cdot \mathbf{x}^t\in \left(\begin{array}{c} m_1\No^*\\ \vdots \\ m_n\No^* \end{array}\right)  \end{equation}
We will use this terminology to describe our example.

 Let $R$ be a (noetherian) semilocal ring such that  $M=\dim _{\varphi} (V^*(R))$ is  the set of solutions of the system given by the two congruences
\[2x+3y\in 5\No^*;\quad x+2y\in 2\No^*\]
Since $(10,0)$ and $(0,5)\in M$, $R$ is semi-semiperfect. Then $V(R)$ is generated by $G=\{(2,2),(6,1),(10,0),(0,5)\}$, and  $V^*(R)$ is generated by $G\bigcup \{(\infty ,0), (0,\infty), (1,\infty)\}$. Notice that a projective right $R$-module $P$ such that $\dim _{\varphi} (\langle P \rangle )=(1,\infty)$ is not a direct sum of finitely generated  modules but so is $P^2$. 
\end{Ex}

\section{Fair sized projective modules} \label{sectionfair}

 Let $R$ be a ring, and  let $M$ be a right $R$-module.
Consider the set
\[\mathcal{I}(M)=\{I\unlhd R\mid M/MI\mbox{ is finitely generated}\}\]
Note that
 $R\in \mathcal{I}(M)$, and that
 $\{0\}\in \mathcal{I}(M)$ if and only if $M$ is finitely generated.

In general, if $I\in \mathcal{I}(M)$ and $J$ is a two-sided ideal of
$R$ such that $I\subseteq J$ then $J\in \mathcal{I}(M)$.

When $M$ is a countably generated projective module there is a way to
produce elements of $\mathcal{I}(M)$ using an idempotent matrix
defining $M$.

\begin{Lemma}\label{descending} Let $R$ be a ring let $A=(a_{i,j})$
be a countable column finite matrix with entries in $R$ such that
$A^2=A$, and set $P=AR^{(\N)}$. For any  $k\ge 0$, let $I_k=\sum
_{i>k, j\in \N}a_{ij}R$. Then,
\begin{itemize}
\item[(i)] $\mathrm{Tr} \, (P)=RI_0$.
\item[(ii)] For any $k\ge 0$ and any $0<i\le k$, $$\sum _{j \in
\N}a_{ij}R+RI_k/RI_k=\sum _{1\le \ell\le k}a_{i\ell}R+RI_k/RI_k.$$
In particular if $0\le n\le k$, $I_n+RI_k/RI_k$ is finitely
generated.
\item[(iii)] For any $k\ge 0$, $RI_k\in \mathcal{I} (P)$.
\item[(iv)] For any $k\ge 0$, there exists $n_k>k$ such that
$I_{n_k}I_k=I_{n_k}$.
\end{itemize}
\end{Lemma}

\begin{Proof} Statement  $(i)$ is clear. To prove  $(ii)$ note  that, since $A$ is an idempotent
matrix, for any $j\in \N$
\[a_{ij}=\sum _{\ell \in \N} a_{i\ell}a_{\ell j} \in \sum _{\ell\le
k}a_{i\ell}a_{\ell j}+RI_k\subseteq \sum _{1\le \ell\le
k}a_{i\ell}R+RI_k/RI_k.\]

$(iii)$. The projective module $P$ is generated by the columns of
$A$. The $R/RI_k$ projective module $P/PI_k$ is isomorphic to the
module generated by the columns of the matrix $(a_{ij}+RI_k)$ which
is an idempotent matrix with entries in $R/RI_k$ such that only the
first $k$ rows may be different from zero. So $P/PI_k$ is a direct
summand of $\left(R/RI_k\right)^k$, therefore it is a finitely
generated projective module.

$(iv)$. Since the matrix $A$ is column finite, for any given $k$
there exists $n_k>k$ such that $a_{ij}=0$ for $i\ge n_k$ and $j\in
\{1,\dots ,k\}$.

Now if $a_{ij}$ is such that $i> n_k$ then either $a_{ij}=0$ if
$j\in \{1,\dots ,k\}$ or
\[a_{ij}=\sum _{\ell >k} a_{i\ell}a_{\ell j}\in I_{n_k} I_k.\]
This implies that $I_{n_k} I_k=I_{n_k}$.
\end{Proof}

\begin{Lemma}\label{star} \cite[Lemma~2.4]{fair} Let $R$ be a ring. Let $P$ be a countably generated projective right
$R$-module. Let $A=(a_{ij})$ be a countable column finite idempotent
matrix such that $P\cong AR^{(\N)}$. For any $k\in \N$ set $I_k=\sum
_{j\in \N \, i>k}a_{ij}R$. Then
\begin{itemize}
\item[(i)] For any $I\in \mathcal{I}(P)$, there exists $k\in \N$
such that $RI_k\subseteq I$.
\item[(ii)] $\mathcal{I}(P)$ is closed under finite intersections.
\item[(iii)] If $\mathcal{I}(P)$ has minimal elements then it has a
unique minimal element.
\item[(iv)]  $\mathcal{I}(P)$ has minimal elements  if and only if
there exists $k_0\in \N$ such that, for any $\ell \in \N$,
$RI_{k_0}=RI_{k_0+\ell}$. In this case $I_{k_0}=I_{k_0}^2$ and
$RI_{k_0}$ is the minimal element of $\mathcal{I}(P)$.
\end{itemize}
\end{Lemma}

\begin{Proof} $(i)$. Since  $AR^{(\N)}$ is generated by the columns $C_1,\dots, C_n,\dots $ of
the matrix $A$,  if $\overline{P}=AR^{(\N)}/AI^{(\N)}$ is finitely
generated there exists $n$ such that $S=\{C_1+ AI^{(\N)}, \dots
,C_n+ AI^{(\N)}\}$ generate $\overline{P}$.  As $A$ is column
finite, there exists $k$
 such
that, $a_{ij}=0$ for $j\in \{1,\dots ,n\}$ and $i>k$. Since for any
$\ell
>n$, $C_\ell +AI^{(\N)}$ is an $R$-linear combination of the elements in
$S$, we deduce that $a_{ij}\in I$ for any $j\in \N$ and any $i>k$.
So that $RI_k\subseteq I$.

Statement $(ii)$ follows from $(i)$ and from the fact that, for any
$k$, $RI_k\in \mathcal{I} (P)$. Statement $(iii)$ follows from
$(ii)$, and $(iv)$ is a consequence of $(ii)$ and
Lemma~\ref{descending}$(iv)$.
\end{Proof}

Following the idea of \cite[proof of Proposition
4.2]{emmanouiltalelli} we can prove further closure properties of the set $\mathcal{I} (P)$.

\begin{Lemma} \label{criteriafg} Let $R$ be a ring and let $P$ be a right $R$-module and let $P_0$ be a submodule of $P$.
Let $I$ be a two-sided ideal of $R$. Then
\begin{itemize}
\item[(i)] $P=P_0+PI$ if and only if, for any right $R/I$-module $M$, $\mathrm{Hom}_R(P/P_0,M)=0$
\item[(ii)] Assume $P$ is projective and $P_0$ is finitely generated
then if $\mathrm{Hom}_R (P/P_0,R/I)=0$ then $P/PI$ is a finitely
generated projective right $R$-module.
\end{itemize}
\end{Lemma}

\begin{Proof} To prove both statements we have to keep in mind the
following particular case of the
$\mathrm{Hom}$-$\otimes$-adjunction: for any right $R/I$-module $M$
\[\mathrm{Hom}_R(P/P_0,M)\cong
\mathrm{Hom}_R(P/P_0,\mathrm{Hom}_{R/I}(R/I,M))\cong\]\[\cong
\mathrm{Hom}_{R/I}(P/P_0\otimes _R
R/I,M)\cong\mathrm{Hom}_{R/I}(P/(P_0+PI),M)\]

Then $(i)$ is clear. To prove $(ii)$, notice that the assumption and the above formula 
imply that $$\mathrm{Hom}_{R/I}(P/(P_0+PI),R/I)=0.$$ Notice that 
$P/(P_0+PI)\cong \left(P/PI\right)/\left(P_0+PI/PI\right)$ and that
$P_0+PI/PI$ is a finitely generated submodule of the projective
$R/I$-module $P/PI$ . Now we follow the argument of \cite[Lemma 4.1]{emmanouiltalelli} to conclude. 

There is a set $X$ an a splitting embedding $f\colon P/PI\to (R/I)^{(X)}$; since $P_0+PI/PI$ is finitely generated, there exists a finite subset $Y$ of $X$ such that $f(P_0+PI/PI)\subseteq (R/I)^{Y}$. For any $x\in X$, let $p_x\colon (R/I)^{(X)}\to R/I$ denote the projection onto the $x$-component. The hypothesis implies that $p_xf=0$ for any $x\in X\setminus Y$, hence $P/PI$ is isomorphic to a direct summand of $(R/I)^{Y}$ and, hence, it is finitely generated.
\end{Proof}

\begin{Prop} Let $R$ be a ring and let $P$ be a projective right
$R$-module. Fix a finitely generated submodule $P_0$ of $P$ and set
\[\mathcal{I} (P_0,P)=\{I \mbox{  ideal of  } R\mid P=P_0+PI\}.\]
Then
\begin{itemize}
\item[(1)] $\mathcal{I} (P_0,P)$ is closed under  products.
\item[(2)] if $\{I_i\}_{i\in \Lambda}$ is a family of ideals in
$\mathcal{I} (P_0,P)$ then $\bigcap _{i\in \Lambda} I_i\in
\mathcal{I}(P)$.
\item[(3)] If $I\in \mathcal{I}\, (P)$ then $I^\omega =\bigcap
_{n\ge 1}I^n\in \mathcal{I}\, (P)$.
\end{itemize}
\end{Prop}

\begin{Proof} $(1)$. Let $I$, $J$ be ideals in $\mathcal{I}
(P_0,P)$. In view of Lemma \ref{criteriafg}(i), we must show that
$\mathrm{Hom}_R(P/P_0,M)=0$ for any right $R/IJ$ module $M$.
To this aim we apply the functor $\mathrm{Hom}_R(P/P_0,-) $ to the
short exact sequence
\[0\to MI\to M\to M/MI\to 0\] to obtain the exact sequence
\[0\to \mathrm{Hom}_R(P/P_0,MI)\to \mathrm{Hom}_R(P/P_0,M)\to \mathrm{Hom}_R(P/P_0,M/MI).\]
Since $MI$ is a right $R/J$-module, $\mathrm{Hom}_R(P/P_0,MI)=0$
and, since $M/MI$ is a right $R/I$-module,
$\mathrm{Hom}_R(P/P_0,M/MI)=0$. Hence $\mathrm{Hom}_R(P/P_0,M)=0$ as
wanted.

$(2)$. Applying the functor $\mathrm{Hom}_R(P/P_0,-) $ to the exact
sequence \[0\to R/\bigcap _{i\in \Lambda}I_i\to \prod_{i\in
\Lambda}R/I_i\] we obtain
\[0\to \mathrm{Hom}_R(P/P_0,R/\bigcap _{i\in \Lambda}I_i)\to \mathrm{Hom}_R(P/P_0,\prod_{i\in
\Lambda}R/I_i)\cong \prod_{i\in
\Lambda}\mathrm{Hom}_R(P/P_0,R/I_i)=0.\] Therefore,
$\mathrm{Hom}_R(P/P_0,R/\bigcap _{i\in \Lambda}I_i)=0$ and then the
statement follows from Lemma~\ref{criteriafg}(ii).

Statement $(3)$ follows from $(1)$ and $(2)$.
\end{Proof}

\begin{Def} \emph{(\cite{fair})} Let $R$ be a ring.  A countably generated projective
right $R$-module $P$ is \emph{fair sized} if $\mathcal{I}\, (P)$ has
a minimal element. \end{Def}

\begin{Prop} \label{minimal} Let $R$ be a semilocal ring, and let $P$ be a countably
generated projective right $R$-module then $P$ is fair sized.
\end{Prop}

\begin{Proof} We can assume that $P$ is not finitely generated and $\mathcal{I}\,
(P)\neq \{R\}$.

Let $\mathcal{S}=\{\mathcal{M}_1,\dots ,\mathcal{M}_k\}$ denote the
set of all maximal two-sided ideals of $R$. If $R\neq I\in
\mathcal{I}\, (P)$, then there exists $i\in \{1,\dots ,k\}$ such
that $I\subseteq \mathcal{M}_i$. In particular, $\mathcal{M}_i\in
\mathcal{I}\, (P)$.

Reindexing the elements of $\mathcal{S}$,  if necessary, we may
assume that there exists $r\ge 1$ such that $P/P\mathcal{M}_i$ is
finitely generated if $i\le r$ and $P/P\mathcal{M}_i$ is infinitely
generated for $i>r$. By Lemma~\ref{star}(ii),
$J=\mathcal{M}_1\bigcap \cdots \bigcap\mathcal{M}_r\in \mathcal{I}
(P)$.

Fix $A=(a_{ij})$ a column finite idempotent matrix such that $P\cong
AR^{(\N)}$. For $k>0$, set $I_k=\sum _{i<k\, j\in \N} a_{ij}R$. By
Lemma~\ref{star}(i), there exists $k_0 $ such that for any $\ell\ge
k_0$, $RI_\ell \subseteq J$. Moreover, in view of
Lemma~\ref{descending} and the above remarks, for any $\ell\ge k_0$,
$J(R/RI_\ell)=J/RI_\ell$.

By Lemma~\ref{descending}, there exists $k_0'>k_0$ such that $I_{k_0
'}I_{k_0} =I_{k_0 '}$. Let $\ell > k_0 '$. Then
\[M=I_{k_0 '}+RI_{\ell}/RI_{\ell}=\left(I_{k_0
'}+RI_{\ell}/RI_{\ell}\right)\left(I_{k_0
}+RI_{\ell}/RI_{\ell}\right)=MJ(R/RI_{\ell}).\] By
Lemma~\ref{descending}, $M$ is a finitely generated right
$R/RI_{\ell}$-module. By Nakayama's Lemma,  $M=0$ or, equivalently,
$I_{k_0 '}\subseteq RI_{\ell}$. In view of Lemma~\ref{star}, we
deduce that $\mathcal{I}\, (P)$ has a least element as we wanted to
prove.
\end{Proof}

The first part of the next result is due to Sakhaev
\cite[Theorem~4]{Sa}. We give an alternative proof following the
same ideas  used in the proof of Proposition~\ref{minimal}.

\begin{Prop}\label{idempotent}  Let $R$ be a ring. Let $P$ be a countably generated
projective right $R$-module such that $P/PJ(R)$ is finitely
generated. Then $\mathcal{I}\, (P)$ has a minimal element.

In particular, if $P$ is not finitely generated then $J(R)$ contains
a non-zero idempotent ideal.
\end{Prop}

\begin{Proof} We may assume that $P$ is not finitely generated. Fix an idempotent column
finite matrix $A=(a_{ij})$ such that $P\cong AR^{(\N)}$. For any
$k\in \N$ set $I_k=\sum _{j\in \N \, i<k}a_{ij}R$. Since $J(R)\in
\mathcal{I}\, (P)$, by Lemma~\ref{star}(i), there exists $k_0$ such
that $RI_{k_0}\subseteq J(R)$. Now the statement follows using the
same argument as in the proof of Proposition~\ref{minimal}.

 By
Lemma~\ref{star}, there exists  $\ell\in \N$ such that the minimal
element in $\mathcal{I}(P)$ is $RI_{\ell}\neq 0$ which, in addition,
is an idempotent ideal. By the above argument $RI_{\ell}\subseteq
J(R)$.
\end{Proof}

Since a projective module with a semilocal endomorphism ring is
finitely generated modulo its Jacobson radical we deduce the
following corollary.

\begin{Cor} Any projective module with a semilocal endomorphism ring
is fair sized.
\end{Cor}

\begin{Remark} The first example of a countably generated projective module, not
finitely generated, but finitely generated modulo its Jacobson
radical was given in \cite{GS}. Such example is, in fact, over a
semilocal ring.

In \cite{DPP}, a  detailed study of Gerasimov and Sakhaev example
was made. In particular, it  follows from \cite[Theorem~6.8]{DPP},
that such a ring $R$ has an idempotent ideal that it is  trace ideal
of a projective right $R$-module but it is not a trace ideal of a
projective left $R$-module.

In \cite[Question 5.6, Question 6.6]{DPP} it was asked whether
Gerasimov and Sakhaev's example satisfies that $\bigcap_{n\ge 1}
J(R)^n=\{0\}$ and whether such a ring has idempotent ideals
contained in $J(R)$. Propositions~\ref{minimal} and \ref{idempotent}
 answer  both questions in the negative. Notice that \cite[Theorem
 C]{emmanouiltalelli} gives also a negative answer to the first
 question.
 \end{Remark}

Let $R$ be a noetherian ring. Let $P_1$ and $P_2$ be   countably generated, fair sized, projective right $R$-modules such that, for $i=1,2$,  $I$ is the minimal ideal such that $P_i/P_iI$ is finitely generated. Then $P_1\cong P_2$ is and only if $P_1/P_1I\cong P_2/P_2I$ \cite{fair}. This result was the crucial tool to classify projective modules over some classes of noetherian rings: generalized Weyl algebras, semilocal noetherian rings, some lattices. However, as the following example shows, this is no longer true in general.

\begin{Ex}   Let $R$ be a ring such that there is an isomorphism of monoids $\varphi \colon V(R)\to  \mathbb{R}_0$ and such that $\varphi (\langle R\rangle)>\sum _{i\ge 1}\frac 1{i^2}$. The existence of such ring is ensured by a theorem due to Bergman and Dicks \cite{BergmanDicks}.

For any $i\ge 1$, let $P_i$ be a finitely generated projective right $R$-module such that $\varphi (\langle P\rangle)=\frac 1{i^2}$. Then $R^{(\omega)}$ and $P=\oplus _{i\ge 1}P_i$ are uniformly big projective modules. That is, both are fair sized and the minimal ideal such that modulo it they are finitely generated is $I=R$. Hence $0=R^{(\omega)}/R^{(\omega)}I=P/PI$ but $ P$ has no nonzero free direct summands. In particular $P\not \cong R^{(\omega)}$.
\end{Ex}

\section{A particular way to think the reconstruction of the projective module} \label{sectionconstruction}

Let  $T$ be a countable rooted tree with no leafs.  Let $V$ be the
set of vertices of $T$,
 let $v_0\in V$ be the root and let $E$ be the set of vertices. Let $s\colon E\to V$ be the source map and $t\colon E\to V$ be the end
 map. That is, an edge $e$ starts in $s(e)$ and terminates in
 $t(e)$. We assume that for each vertex $v$ there is only a finite number of edges $e$ such that $s(e)=v$.
 
 Consider the distance function $d \colon V \times V \to \N_0$  where, by definition,
$d(v,v')$ is the length of the shortest (unoriented) path from $v$
to $v'$.  For every $i \in \N_0$, set  $V_i = \{v \in V \mid
d(v,v_0) = i\}$, and let $E_i$ be the set of edges connecting a
vertex of $V_i$ and a vertex of $V_{i+1}$.

 Let $R$ be a ring, and  fix a representation of $T$ in $\Mod R$. The representation assigns to each
$v \in V$ a module $M_v$ and to every $e \in E$ a homomorphism $f_e \colon M_{s(e)} \to M_{t(e)}$. Set $M_i = \oplus_{v \in V_i} M_v$. For every
$v \in V_i$ let $\iota^i_v \colon M_v  \to M_i$ be the canonical
embedding and let $\pi ^i_v \colon M_i \to M_v$ be the canonical
projection. For every $i \in \N_0$ there is a unique homomorphism
$f_i \colon M_i \to M_{i+1}$ such that, for every $e \in E_i$,
$\pi ^{i+1}_{t(e)} f_i \iota^i_{s(e)} = f_e$. Therefore we get a direct system
\[M_0 \stackrel{f_0}\to M_1 \stackrel{f_1}\to\cdots
M_k\stackrel{f_k}\to M_{k+1}\cdots \] whose limit fits into
a pure exact sequence
\begin{equation} \label{telescope} 0\to \oplus_{v \in V} M_v \to \oplus_{v \in V} M_v\to \varinjlim
M_i\to 0\end{equation}

Now we introduce some further notation in order to be able to give a  local criteria to ensure that this exact sequence splits.

  For every $v \in V$ let
us denote $S(v)$ the set of successors of $v$, that is $S(v) = \{v'\in V \mid
\exists e \in E\mbox{ such that } s(e) = v, t(e)=v'\}$.  Let $S_2(v) = \bigcup_{v'
\in S(v)} S(v')$. For any vertex $v$ let $M_v' = \oplus_{u \in S(v)} M_{u}$ and let $M_v'' = \oplus_{u \in S_2(v)} M_u$.
For every $u \in S(v)$ and $w \in S_2(v)$ let
$\pi_{u}' \colon M_v' \to M_{u}$ and  $\pi_{w}'' \colon M_{v}'' \to M_w$
be the canonical projections and let $\iota_{u}' \colon M_u \to M_v'$ be the canonical
embedding. Then there exist a
unique homomorphism $g_v \colon M_{v} \to M_{v}'$ 
such that $\pi_{t(e)}' g_v = f_e$  for any  $e \in E$ such that $v = s(e)$. There is also a unique homomorphism $g_v' \colon M_{v}' \to M_{v}''$ such that 
$\pi_{t(e')}'' g_v' \iota_{s(e')}' = f_{e'}$ for any $e' \in E$ satisfying that $s(e') \in S(v)$.

\begin{Lemma}\label{splitm} With the notation above. If,  for every $v \in V$, there exists a homomorphism $h_v\colon
M_v'' \to M_v'$ such that $h_vg_v'g_v = g_v$ then the  
sequence (\ref{telescope}) splits.
\end{Lemma}

\begin{Proof}  Notice that, for every $i \in \N_0$,   $M_{i+1} = \oplus_{v \in V_i}
M_v'$ and $M_{i+2} = \oplus_{v \in V_i} M_{v}''$. Therefore,  for every $i \in
\N_0$, there is
a homomorphism $g_{i+1} \colon M_{i+2} \to M_{i+1}$ given by $g_{i+1} =
\oplus_{v \in V_i} h_v$. Using the condition $h_vg_v'g_v = g_v$ it
is easy to see that $g_{i+1}f_{i+1}f_i = f_i$ holds. By Lemma \ref{split}, the sequence (\ref{telescope}) splits.
\end{Proof}

\begin{Ex}\label{whiteheadgraph}
 Suppose that $I$ is a nonzero idempotent ideal of $R$ generated as a left ideal by $a_1,a_2,\dots,a_\ell$. Then a construction of a projective
module having the trace ideal $I$ can be explained as follows: Let
$T$ be a rooted tree such that, for every vertex $v$, $\vert S(v)\vert
= n$.  For every  $v \in V$, set $M_v = R$ and, if $e_1,\dots,e_\ell
\in E$ are the edges such that $s(e) = v$, set $f_{e_i} \colon R \to
R$ to be left multiplication by $a_i$. So that, the representation can be visualized as a repetition at each vertex of the following basic tree of height $1$:

$$
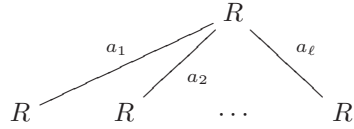
\xymatrix{ & &R\ar@{-}[dr]^{a_\ell} \ar@{-}[dll] _{a_1}\ar@{-}[dl] ^{a_2} &\\ R& R &\cdots &R}$$
so that, for each $v\in V$, $g_v\colon R\to R^\ell$ is given by $g_v(r)=\begin{pmatrix}a_1\\ \vdots \\ a_\ell \end{pmatrix}r$ for any $r\in R$. Set $A=\begin{pmatrix}a_1\\ \vdots \\ a_\ell \end{pmatrix}$ and $B=\left( \begin{array}{ccc} A&\dots &0\\
\vdots & \ddots ^{\ell)} &\vdots \\ 0& \dots & A \end{array}\right)$. Then $g'_v\colon R^\ell \to (R^\ell)^\ell$ is given by $g'_v(r_1,\dots ,r_\ell)=B \begin{pmatrix}r_1\\ \vdots \\ r_\ell \end{pmatrix}$. By Lemma \ref{independent}, there exists $C\in M_{\ell \times \ell ^2}(R)$ such that $CBA=A$. Now the homomorphism $h_v\colon (R^\ell)^\ell\to R^\ell$ given by left multiplication by $C$  satisfies that $h_vg_v'g_v = g_v$. Therefore,
by Lemma \ref{splitm}, the   limit of the direct system
\[M_0 \stackrel{f_0}\to M_1 \stackrel{f_1}\to\cdots
M_k\stackrel{f_k}\to M_{k+1}\cdots   \] 
where, for any $k\ge 0$, $M_k=R^{\ell ^{k}}$,  is
 a projective module $P$. By Proposition \ref{directsystem}, the trace  of
$P$ is the ideal $I$.

Observe  that, for any $k\ge 0$,
$$P\cong \varinjlim _{i\ge k}M_i\cong (\varinjlim _{i\ge
0}M_i)^{\ell ^k}=P^{\ell ^k}\neq \{ 0\} .$$ In particular, if $R$ is a right noetherian
ring, it follows that such projective module $P$ cannot have finite Goldie dimension and, hence, it is always $I$-big, that is, $P$ is not finitely generated, and if $J$ is a right ideal of $R$
such that $P/PJ$ is finitely generated then $P=PJ$ and, hence, $RJ\supseteq I$.
\end{Ex}

Our next construction will use the following lemma.

\begin{Lemma} \label{independentidempotent} Let $R$ be a ring, and
let $I_1\supseteq I_2\supseteq I_3$ be two-sided ideals of $R$ such
that $I_3$ is idempotent and ${}_RI_3$ is finitely generated. Assume
that $R/I_3$ is semisimple  and, for $i\le j$, let $e_{ij}\in R$ be
such that $e_{ij}+I_j$ is a central idempotent of $R/I_j$ generating
$I_i/I_j$. Then,
\begin{itemize}
\item[(i)] $I_1$ and $I_2$ are also idempotent ideals which are finitely
generated as left $R$-modules.
\item[(ii)] Let $G\in M_{m\times 1}(R)$ and $H\in M_{n\times 1}(R)$ be column matrices
such that their entries generate ${}_RI_2$ and ${}_RI_3$,
respectively. Set $A=\begin{pmatrix}e_{1,2}\\ G\end{pmatrix}$, and
\[B=\left(\begin{array}{cccc}Y&&\dots &0\\0& Z&\dots &0\\&&\ddots ^{m)} &\\
&&&Z\end{array}\right)\] where $Y=\begin{pmatrix}e_{1,3}\\
H\end{pmatrix}$ and $Z=\begin{pmatrix}e_{2,3}\\
H\end{pmatrix}$. Then there exists a matrix $C$ of suitable size and
with entries in $R$ such that $CBA=A$.
\end{itemize}
\end{Lemma}

\begin{Proof} $(i)$. This statement follows easily from the fact
that ${}_RI_3$ is finitely generated and $R/I_3$ is semisimple
artinian.

$(ii)$. Since $e_{1,2}R\subseteq R e_{1,2}+I_2$ and $I_2=\sum _{i=1}^mRg_i$, where
$g_1,\dots, g_m$ denote the entries of $G$, it follows as in Example~\ref{whiteheadgraph} that there exists $D\in M_{(m+1)\times (m+1)}(R)$ such that $DA=A$.  Since $e^2_{1,2}-e_{1,2}\in I_2$,  $D$ can be taken of the form
$D=(A\vert D')$ where $D'$ has entries in $I_2$.

Since the entries in $Y$ generate $I_1$ and the entries in $Z$
generate $I_2$, one can  proceed as in Lemma \ref{independent} to
show that there exists a suitable matrix $C$ such that  $D=CB$.
\end{Proof}

We are also interested in the following variation of Lemma
\ref{independentidempotent}.

\begin{Rem} \label{variation} In the same situation as in Lemma
\ref{independentidempotent}, fix $c $, $d\in \mathbb{N}$. Let
$A_c=\begin{pmatrix}e_{1,2}\\ \vdots \, {\scriptstyle{c)}}\\e_{1,2} \\
G\end{pmatrix}$. Let \[B_{cd}=\left(\begin{array}{cccccc}Y_c&&&&\dots &0\\&\ddots &&&&\\&&Y_c&&&\\0&&& Z_d&\dots &0\\&&&&\ddots ^{m)} &\\
&&&&&Z_d\end{array}\right)\] where $Y_c=\begin{pmatrix}e_{1,3}\\ \vdots \,{\scriptstyle{c)}}\\e_{1,3}\\
H\end{pmatrix}$ and $Z_d=\begin{pmatrix}e_{2,3}\\ \vdots \, {\scriptstyle{d)}}\\e_{2,3}\\
H\end{pmatrix}$. Then there exists a matrix $C_{cd}$ of suitable
size and with entries in $R$ such that $C_{cd}B_{cd}A_c=A_c$.
\end{Rem}

From now on we assume that $R$ is a   ring having a descending chain
of ideals $I_0 \supseteq I_1 \supseteq \cdots $ such that
\begin{enumerate}
\item[1)] For every $i \in \N_0$ the ring $R/I_i$ is semisimple artinian.
\item[2)] For every $i \in \N_0$ the ideal $I_i$ is idempotent and ${}_RI_i$ is finitely generated.
\end{enumerate}

For every $j \in \N$ let $G_j$ be a finite set of generators
${}_RI_j$, and for every $i,j \in \N_0, i < j$ let $e_{i,j}$ be
chosen such that $e_{i,j} + I_j$ is a central idempotent of $R/I_j$
generating $I_i/I_j$. Further let $m_i = |G_i|$ for every $i \in
\N$.

Now we construct a tree $T$ and a representation of $T$ by
$R$-modules. We proceed by   induction on the distance from the
root and we construct $T$ and the representation simultaneously. We
keep the notation from the beginning of the section. For every $v
\in V$ we put $M_v = R$. The root $v_0$ of $T$ has $1+m_1$
successors. If $e_0,\dots,e_{m_1}$ are the edges adjacent to $v_0$
then $f_{e_0}$ is given by left multiplication by $e_{0,1}$, and $f_{e_1},\dots,f_{e_{m_1}}$ are given, respectively, by left multiplication
by the $m_1$ elements of $G_1$. That is, if $G_1=\{g_{1, 1},\cdots , g_{1, m_1}\}$ the result of this first step can be visualized as 
$$
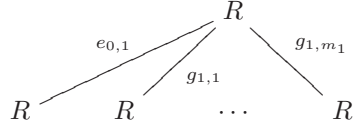
\xymatrix{ & &R\ar@{-}[dr]^{g_{1, m_1}} \ar@{-}[dll] _{e_{0,1}}\ar@{-}[dl] ^{g_{1, 1}} &\\ R& R &\cdots &R}$$

Now suppose we have constructed the subtree and the representation
on the subtree $T'$ given by the set of vertices $\bigcup_{0\leq j \leq
i} V_{j}$ for some $i \in \N$ in such a way that for  $v \in V_i$:
\begin{enumerate}
\item[(i)]   
there is a unique edge $e_v$ in $T'$ such that $v = t(e_v)$;
\item[(ii)] $f_{e_v}$ is given by
\begin{enumerate}
\item[(ii.a)] left  multiplication by $e_{j,i}$   for some $j<i$ or
\item[(ii.b)]  by left multiplication by an
element from $G_{i}$. 
\end{enumerate}
\end{enumerate}
Now for $v \in V_i$, we add exactly $1+m_{i+1}$ successors of $v$
to $T'$ and the corresponding edges say $e_0,\dots,e_{m_{i+1}}$. Then we
set $f_{e_0}$ to be given by left multiplication by $e_{j,i+1}$ in the situation $(ii.a)$, otherwise we
set $f(e_0)$ to be  left multiplication by $e_{i,i+1}$. Further, we set 
$f_{e_1},\dots,f_{e_{m_{i+1}}}$ to be given by  left multiplication
by the elements of $G_{i+1}$. Repeating this procedure for every $v
\in V_i$ completes the inductive step.

By Lemma \ref{independentidempotent}, for every $v \in V$ there
exists $h_v \colon M_{v}'' \to M_{v}'$ such that $h_vg_v'g_v = g_v$. Hence, by Lemma~\ref{splitm}, the  direct limit of the system $M_0 \stackrel{f_0}\to M_1
\stackrel{f_1}\to\cdots M_k\stackrel{f_k}\to M_{k+1}\cdots$ given by
the representation of $T$ is a countably generated projective right module that we denote by $P$.

The module  $P$ satisfies
the following properties.

\begin{Lemma} \label{invariantsofp}
\begin{itemize}
\item[(i)] The trace ideal of $P$ is $I_0$.
\item[(ii)] For every $k\ge 0$, $$P\otimes _R R/I_{k+1}\cong \left(I_0/I_{k+1}\right)\oplus \left( I_1/I_{k+1}\right)^{m_1}\oplus \left( I_2/I_{k+1}\right)^{(m_1+1)m_2}\oplus
\cdots \oplus \left( I_k/I_{k+1}\right)^{(m_1+1)(m_2+1)\dots m_k}
.$$ In particular, $P/PI_k$ is finitely generated for any $k\ge 0$.
\item[(iii)] For every $k\ge 0$,
$P\cong P_{k0}\oplus \cdots \oplus P_{kk}$,  where $\mathrm{Tr}\,
(P_{ki})=I_i$ for $i=0,\dots ,k$.
\item[(iv)] If $R_R$ is Noetherian and, for any $k\ge 0$, $I_k\varsupsetneq
I_{k+1}$ then $P$ is not fair sized.
\end{itemize}
\end{Lemma}

\begin{Proof} Statement $(i)$ is clear from the construction and
Proposition \ref{directsystem} because the entries of the transition
maps of the direct limit generate the ideal $I_0$.

To prove statement $(ii)$, note  that 
$$P/PI_{k+1}\cong P\otimes _R R/I_{k+1}\cong \varinjlim _{i\ge k} M_i \otimes _R R/I_{k+1}\cong \varinjlim _{i\ge k} M_i/M_iI_{k+1}.$$
The statement follows counting the number of  semisimple factors in
this direct limit.

Statement $(iii)$ follows in a similar way, taking into account that
the direct system $\{M_i\}_{i\ge k}$ is the direct sum of $k+1$
direct systems, each one with limit a projective module. By
Proposition \ref{directsystem} the traces of such projective modules
are the ideals $I_0,\dots ,I_k$, respectively.

$(iv)$. Assume there is an ideal of $R$ minimal with respect to the
property that $P/PI$ is finitely generated. Since $R/I$ is a right
noetherian ring, $P/PI$ is a module with finite Goldie dimension
$n$, say. In view of $(iii)$, this implies that $P_{nn}/P_{nn}I=0$.
Equivalently, $I\supseteq I_n$. Then statement $(ii)$ contradicts
the minimality of $I$.
\end{Proof}

Making a variation of the construction of $P$ we have the following
extension of Lemma~\ref{invariantsofp}.

\begin{Lemma} \label{invariantsofps} Let $R$ be a ring. For any  sequence $s$ of elements in $\N$
there exists a countably generated projective right $R$-module $P_s$ with
trace $I_0$, and satisfying the following properties:
\begin{itemize}
\item[(i)] For every $k\ge 0$, $$P_s\otimes _R R/I_{k+1}\cong \left(I_0/I_{k+1}\right)^{\alpha_0}
\oplus \left( I_1/I_{k+1}\right)^{\alpha_1}\oplus   \cdots \oplus
\left( I_k/I_{k+1}\right)^{\alpha_k},$$
where $\alpha_0,\dots,\alpha_k \in \N_0$ determine first $k+1$ terms of the sequence s. In
particular, $P_s/P_sI_k$ is finitely generated for any $k\ge 0$.
\item[(ii)] If, for any $k\ge 0$, $I_k\varsupsetneq
I_{k+1}$ then $P_s\cong P_{s'}$ implies $s=s'$.
\item[(iii)] For every $k\ge 0$,
$P_s\cong Q_{k0}\oplus \cdots \oplus Q_{kk}$, where $\mathrm{Tr}\,
(Q_{ki})=I_i$ for $i=0,\dots ,k$.
\item[(iv)] If $R_R$ is Noetherian and, for any $k\ge 0$, $I_k\varsupsetneq
I_{k+1}$ then $P_s$ is not fair-sized.
\end{itemize}
\end{Lemma}

\begin{Proof} Let $c_0,c_1,\dots \in \N$ be the elements of the sequence $s$. Modify the construction before  Lemma~\ref{invariantsofp} as follows: Set $M_v = R$ for every $v \in V$, and let
$a_{0,1},\dots,a_{0,c_0},e_1,\dots,e_{m_1}$ be the edges starting at
$v_0$. Then $f_{a_{0,i}}$  is given by left multiplication by $e_{0,1}$
and $f_{e_1},\dots,f_{e_{m_1}}$ are given, respectively, by  left multiplication
by the $m_1$ elements of $G_1$. For the inductive step, assume $i>0$ and fix $v \in V_i$. Then $v = t(e_v)$ for a unique
$e_v \in E$. Put $j = i$ provided $f_{e_v}$ is given by multiplication by an element in $G_i$ and put $j = k$ if $f_{e_v}$ is given by
left  multiplication by $e_{k,i}$ for some $k<i$. 
If $j < i$ the vertex $v$ has $1+m_{i+1}$ successors, let us denote 
$e_0,e_1,\dots,e_{m_{i+1}}$ the edges starting at $v$. Then $f_{e_0}$
is given by  left multiplication by $e_{j,i+1}$ and 
$f_{e_1},\dots,f_{e_{m_{i+1}}}$
are given by  left multiplication by the different elements of $G_{i+1}$.
If $i = j$ the vertex $v$ has
$c_i+m_{i+1}$ successors, let us denote
$a_{i,1},\dots,a_{i,c_i},e_1,\dots,e_{m_{i+1}}$ the edges starting
at $v$. Then $f_{a_{i,1}},\dots,f_{a_{i,c_i}}$ are given by  left
multiplication by $e_{i,i+1}$ and $f_{e_1},\dots,f_{e_{m_{i+1}}}$
are given by  left multiplication by the elements from $G_{i+1}$.

Let $P_s$ be the direct limit of the telescope constructed from this
representation. Using Remark~\ref{variation} and Lemma~\ref{splitm}
it follows that $P_s$ is countably generated projective. Proposition~\ref{directsystem} ensures that the trace of $P_s$ is $I_0$.

Statement $(i)$ follows as Lemma~\ref{invariantsofp}$(ii)$. To prove
statement $(ii)$, let $c'_0,c'_1,\dots \in \N$  be the elements of a
sequence $s'$. Notice that if $P_s\cong P_{s'}$ then, for any $k\ge
0$ there is an isomorphism between the finitely generated semisimple
factors $S=(P_s/P_sI_{k+1})I_k$ and $S'=(P_{s'}/P_{s'}I_{k+1})I_k$.
Therefore, by $(i)$, $c_k=c'_k$ for any $k\ge 0$.

Statements $(iii)$ and $(iv)$ are proved as their counterparts in
Lemma~\ref{invariantsofp}.
\end{Proof}

\section{An application to FCR algebras} \label{sectionfcr}

In this section we present a context in which the constructions of
the previous section can be done: the FCR-algebras. They were introduced in
\cite{KS}, and we start giving
a brief overview of their main properties.

Let $k$ be a field. A $k$-algebra $R$ is said to be an
FCR-algebra provided every finite dimensional representation of $R$ is
completely reducible and  the intersection of the kernels of all
finite dimensional representations of $R$ is zero. Basic examples of
FCR-algebras are the universal enveloping algebras of
finite-dimensional semisimple Lie algebras over a field of
characteristic zero and their quantum deformations. For a list of  examples see \cite{KSW}.

To explain better this setting we recall the following result from
\cite{KSW}.

\begin{Th} \emph{\cite[Theorem~3]{KSW}} \label{KSW} Let $R$ be a $k$-algebra where $k$ is an arbitrary field.
Then the following assertions are equivalent:
\begin{itemize}
\item[(i)] Every finite dimensional representation of $R$ is
completely reducible.
\item[(ii)] Every two-sided ideal $I$ of $R$ of finite
codimension is an idempotent ideal.
\item[(iii)] If $\mathcal{M}_1$ and $\mathcal{M}_2$ are two-sided
maximal ideals of finite codimension then
\[\mathcal{M}_1\bigcap \mathcal{M}_2=\mathcal{M}_1\cdot
\mathcal{M}_2=\mathcal{M}_2\cdot \mathcal{M}_1.\]
\end{itemize}
\end{Th}

 It follows from Theorem \ref{KSW}, that every FCR-algebra $R$ has a system of nonzero idempotent ideals
${\mathcal S}$ such that for every $I \in {\mathcal S}$ the ring
$R/I$ is semisimple artinian and $\bigcap _{I\in \mathcal S} I= 0$.

It seems to be quite challenging to control countably generated
projective modules over  a (noetherian) FCR-algebra $R$. There are plenty of semisimple artinian factors of $R$ where we can
try to distinguish non isomorphic projective $R$-modules.
Unfortunately, these semisimple factors do not provide full
information about projective modules over $R$. For example, if $R =
U({\rm sl}_2(\C))$, then there exists a projective right ideal $I$
which is stably free but not free (see \cite{S}). Therefore, $I \otimes_R
R/I_k$ is a free module for every $k \in \N$.

Since a noetherian  FCR-algebra has plenty of idempotent
ideals, it follows from Corollary~\ref{tracefg} that these ideals
are traces of countably generated projective modules. In the case
$R$ is the enveloping algebras of a semisimple Lie algebra, the only
trace ideal of a nonzero finitely generated projective module is
$R$. This is a result proved in Puninski's notes \cite[Lemma
8.1]{puninski} where it is attributed to Stafford. Since
\cite{puninski} has remained unpublished until now, we include a
proof of this result.

\begin{Lemma} \label{fglie} Let $R$ be the enveloping algebra of a finite
dimensional Lie algebra over a field $k$. Then every non-zero
finitely generated projective right $R$-module is a generator.

In particular, the only idempotent ideals of $R$ that are   traces of
 finitely generated projective right modules are $0$ and $R$.
\end{Lemma}

\begin{Proof} By a result of Quillen \cite{macconnellrobson}, every
non-zero finitely generated projective module $P$ is stably free.
Therefore, there exist $n$ and $m$ such that $P\oplus R^n\cong R^m$.
Since $R$ is noetherian, $R_R$ has finite Goldie dimension, hence
$n< m$.

Let $I$ be the trace ideal of $P$. Since $P=PI$, the above
isomorphism, implies that there is an isomorphisms of right
$R$-modules $(R/I)^n\cong (R/I)^m$. Since $R_R$ is noetherian,   $I=R$.
\end{Proof}

\begin{Th}
Let $R$ be an indecomposable  non finite dimensional, noetherian,
$FCR$-algebra over a field $k$. Let $I_0 \subseteq R$ be an ideal of
finite codimension. Then there are uncountably many non-isomorphic
countably generated projective $R$-modules  that are not fair-sized
and with trace ideal contained in $I_0$.
\end{Th}

\begin{Proof} The hypothesis ensure that we can construct a
strictly descending chain
\[I_0\supsetneqq I_1\supsetneqq \cdots \supsetneqq I_k\supsetneqq
\cdots \] of idempotent ideals of $R$. Now the result follows as an
application of Lemma~\ref{invariantsofps}.
\end{Proof}

\end{document}